%% file: pavarxive.tex
\documentclass[12pt]{article}%
\usepackage{makeidx}
\usepackage{amsmath}
\usepackage{float}
\usepackage{graphicx}
\usepackage{amscd}
\usepackage{amsthm}
\usepackage{amsxtra}
\usepackage{verbatim}
\usepackage[size=2em,Postscript=dvips]{diagrams}
\usepackage{amsfonts}
\usepackage{amssymb}
\usepackage{times}%
\diagramstyle[labelstyle=\scriptstyle]
\theoremstyle{plain}
\newtheorem{X}{X}[section]
\newtheorem{theorem}[X]{Theorem}
\newtheorem{proposition}[X]{Proposition}
\newtheorem{lemma}[X]{Lemma}
\newtheorem{corollary}[X]{Corollary}

\theoremstyle{definition}
\newtheorem{definition}[X]{Definition}
\newtheorem{example}[X]{Example}
\newtheorem{plain}[X]{}

\newtheorem*{note}{Notes}
\newtheorem{remark}[X]{Remark}

\input{def01.tex}

\begin{document}

\title{Periods of Abelian Varieties}
\date{February 27, 2003. Submitted version.}
\author{J.S. Milne\thanks{Partially supported by the National Science Foundation.}}
\maketitle

\begin{abstract}
We prove various characterizations of the period torsor of abelian varieties.

\end{abstract}
\tableofcontents

\renewcommand{\theequation}{\arabic{equation}} \pagestyle{myheadings}
\markboth{INTRODUCTION}{INTRODUCTION} \vspace{0.2in}

\addcontentsline{toc}{section}{Introduction}

\subsection{Introduction}

For an abelian variety $A$ over $\mathbb{Q}{}$, $H^{1}(A(\mathbb{C}%
{}),\mathbb{\mathbb{C}{}})$ has two $\mathbb{Q}{}$-structures, that provided
by singular cohomology $H^{1}(A(\mathbb{C}{}),\mathbb{Q}{})$ and that provided
by de Rham cohomology $\mathbb{H}{}^{1}(A_{\text{Zar}},\Omega_{A/k}^{\bullet
})$. The \emph{periods }of $A$ are the coefficients of the transition matrix
from a basis for one $\mathbb{Q}{}$-space to a basis for the other. It is
known (Deligne 1982) that Hodge classes on $A$ impose algebraic relations on
the periods, and it is conjectured that these are the only such relations.
Thus, there appears to be no hope of obtaining an explicit description of the
periods, but one may still hope to characterize some of the objects attached
to them. The singular and de Rham cohomologies define fibre functors on the
category of motives based on the abelian varieties over $\mathbb{Q}{}$, and
the difference of these functors is measured by the period torsor
$P^{\mathsf{AV}}$. In this paper, we obtain various characterizations of
$P^{\mathsf{AV}}$. Beyond its intrinsic interest, the period torsor controls
the rationality of automorphic vector bundles, and therefore of holomorphic
automorphic forms (see Milne 1988; Milne 1990, III 4).

The first problem one runs into is that $P^{\mathsf{AV}}$ is a torsor for an
affine group scheme $G$ over $\mathbb{Q}{}$ for the flat (more specifically,
the f.p.q.c) topology. Such torsors are classified by the flat cohomology
group $H^{1}(\mathbb{Q}{},G)$ rather than a more familiar Galois cohomology
group. For an algebraic quotient $G_{n}$ of $G$, the two cohomology groups
coincide. In \S 1 we prove that the canonical map%
\[
H^{1}(\mathbb{Q}{},G)\rightarrow\varprojlim H^{1}(\mathbb{Q}{},G_{n})
\]
is surjective, and that the fibre of the map containing the class of a
$G$-torsor $P$ is%
\[
\varprojlim{}^{1}(P\wedge^{G}G_{n})(\mathbb{Q}{})
\]
(nonabelian higher inverse limit). The limits are over the set of algebraic
quotients $G_{n}$ of $G$. We show, moreover, that the full group
$H^{1}(\mathbb{Q}{},G)$ can be interpreted as a Galois cohomology group, but
it is Galois cohomology defined using cochains that are continuous relative to
the inverse limit topology on $\varprojlim G_{n}(\mathbb{Q}{}^{\text{al}})$
(discrete topology on each $G_{n}(\mathbb{Q}{}^{\text{al}})$). It is important
to note that these results depend crucially on the fact that the algebraic
quotients of $G$ form a countable set --- I do not even know how to define a
nonabelian higher inverse limit except for countable coefficient sets. The
remainder of \S 1 reviews results, more-or-less known, concerning nonabelian
higher inverse limits and the classification of morphisms of torsors.

In \S 2, we take up the problem of characterizing $P^{\mathsf{CM}}$, the
period torsor for the category of motives based on abelian varieties of
(potential) CM-type over $\mathbb{Q}{}$. Since, as Deligne has pointed out,
the period torsor $P^{\mathsf{Art}}$ attached to the subcategory of Artin
motives can be explicitly described, it is natural rather to consider the
\textquotedblleft relative\textquotedblright\ problem of characterizing the
morphism $P^{\mathsf{CM}}\rightarrow P^{\mathsf{Art}}$ of torsors. Objects of
this type are classified by the flat cohomology group $H^{1}(\mathbb{Q},{}%
_{f}S)$ where $_{f}S$ is a certain twist of the Serre group. As the Serre
group is commutative, the main result of \S 1 simplifies to an exact sequence%
\[
0\rightarrow\varprojlim{}^{1}{}_{f}S_{n}(\mathbb{Q}{})\rightarrow
H^{1}(\mathbb{Q}{},{}_{f}S)\rightarrow\varprojlim H^{1}(\mathbb{Q}{},{}%
_{f}S_{n})\rightarrow0.
\]
Blasius (unpublished) showed that $\varprojlim H^{1}(\mathbb{Q}{},{}_{f}%
S_{n})$ satisfies a Hasse principle, and Wintenberger (1990) showed that
$\varprojlim H^{1}(\mathbb{Q}{},{}_{f}S_{n})=0$. If $\varprojlim{}^{1}{}%
_{f}S_{n}(\mathbb{Q}{})$ were also zero, then $P^{\mathsf{CM}}\rightarrow
P^{\mathsf{Art}}$ would be characterized up to isomorphism (inducing the
identity on $P^{\mathsf{Art}}$) by its cohomology class. Alas, it is not zero
--- in fact, we show that $\varprojlim{}^{1}{}_{f}S_{n}(\mathbb{Q}{})$ is
uncountable. Our proof of this uses an old theorem of Scholz and Reichardt on
the embedding problem for Galois groups of number fields. It would be
interesting to have more information on $\varprojlim{}^{1}{}_{f}%
S_{n}(\mathbb{Q}{})$.

In \S 3, we take up the problem of chacterizing $P^{\mathsf{AV}}.$ Again, it
is more natural to consider the relative problem of characterizing
$P^{\mathsf{AV}}\rightarrow P^{\mathsf{CM}}$. Among other results, we prove
that the isomorphism class of $P^{\mathsf{AV}}\rightarrow P^{\mathsf{CM}}$ is
uniquely determined by its classes over $\mathbb{Q}{}_{l}$ ($l=2,3,\ldots
,\infty$), about which a great deal is known.

Blasius and Borovoi (1999) study the problem of characterizing $P^{\mathcal{H}%
{}}$ where $\mathcal{H\subsetneqq\,}\AV{}$ is the category of motives based on
a certain class of abelian varieties over $\mathbb{Q}$ whose Mumford-Tate
groups have simply connected derived group. Unfortunately, they make the
(false!) assumption that the flat cohomology groups coincide with the
inverse-limit Galois cohomology groups, i.e., effectively they set all
$\varprojlim^{1}$s equal to zero.\footnote{This alone accounts for most of the
relative brevity of their paper.} If this were their only error, their main
Theorem 1.5 would hold for suitably large finite sets of abelian varieties
satisfying their condition, but they also misidentify the cohomology class
that must be proved trivial\footnote{With their notations, in order to prove
their theorem, they would have to show in \textbf{5.1} that the class of
$\mathcal{P}{}\rightarrow\mathcal{P}{}^{\text{CM }}$in $H^{1}(\mathbb{Q}%
{},(\mathcal{G}{}_{\text{DR}}^{\circ})^{\text{der}})$ is trivial (see (1.10)
below). In fact, they prove only the weaker statement that the image of the
class in $H^{1}(\mathbb{Q}{},\mathcal{G}{}_{\text{DR}}^{\circ})$ is trivial,
which, in fact, is all their hypotheses imply, even when one ignores the
$\varprojlim^{1}$ terms.}.  Nevertheless, by combining (1.14) of this paper
with their arguments, one obtains the following theorem (3.25): the class of
$P^{\mathcal{H}{}}\rightarrow P^{\mathsf{CM}}$ is uniquely determined by its
class over $\mathbb{R}{}$. This observation began my work on this paper.

\subsubsection{Notations and conventions}

\textquotedblleft Variety\textquotedblright\ means geometrically reduced
scheme of finite type over a field. Semisimple algebraic groups are connected
and \textquotedblleft simple\textquotedblright\ for an algebraic group means
\textquotedblleft noncommutative and having no proper closed connected normal
subgroup $\neq1$\textquotedblright. The identity component of a group scheme
$G$ over a field is denoted by $G^{\circ}$. For a connected (pro)reductive
group $G$ over a field, $ZG$ is the centre of $G$, $G^{\text{ad}}$ is the
adjoint group $G/ZG$ of $G$, $G^{\text{der}}$ is the derived group of $G$, and
$G^{\text{ab}}$ is the largest commutative quotient $G/G^{\text{der}}$ of $G$.
The universal covering of a semisimple group $G$ is denoted $\tilde
{G}\rightarrow G$.

The algebraic closure of $\mathbb{Q}{}$ in $\mathbb{C}{}$ is denoted by
$\mathbb{Q}{}^{\text{al}}$, and (except in \S 1) $\Gamma=\Gal(\mathbb{Q}%
{}^{\text{al}}/\mathbb{Q})$. We set $\Gal(\mathbb{C}{}/\mathbb{R}{}%
)=\{1,\iota\}$. A CM-field is any field $E$ algebraic over $\mathbb{Q}{}$
admitting a nontrivial involution $\iota_{E}$ such that $\iota\circ\rho
=\rho\circ\iota_{E}$ for all $\rho\colon E\rightarrow\mathbb{C}{}$.

All categories of motives will be defined using absolute Hodge classes as the
correspondences (Deligne 1982a; Deligne and Milne 1982, \S 6).

We sometimes use $[x]$ to denote an equivalence or isomorphism class
containing $x$. The notation $X\approx Y$ means that $X$ and $Y$ are
isomorphic, and $X\cong Y$ means that $X$ and $Y$ are canonically isomorphic
(or that a particular isomorphism is given).

\newpage\pagestyle{headings}

\section{Preliminaries}

\subsection{Inverse limits}

For an inverse system of groups indexed by $(\mathbb{N},\leq)$,%
\[
(A_{n},u_{n})_{n\in\mathbb{N}{}}=(A_{0}\leftarrow\cdots\leftarrow
A_{n-1}\overset{u_{n}}{\leftarrow}A_{n}\leftarrow\cdots),
\]
define $\varprojlim^{1}A_{n}$ to be the set of orbits for the left action of
the group $\prod_{n}A_{n}$ on the set $\prod_{n}A_{n}$,
\[%
\begin{array}
[c]{ccccc}%
\prod_{n}A_{n} & \times & \prod_{n}A_{n} & \rightarrow & \prod_{n}A_{n}\\
(\ldots,a_{n},\ldots) &  & (\ldots,x_{n},\ldots) & \mapsto & (\ldots
,a_{n}\cdot x_{n}\cdot u_{n+1}(a_{n+1})^{-1},\ldots)\text{.}%
\end{array}
\]
This is a set, pointed by the orbit of $\underline{1}=(1,1,\ldots)$. Note
that
\[
\varprojlim A_{n}=\{a\in%
%TCIMACRO{\tprod }%
%BeginExpansion
{\textstyle\prod}
%EndExpansion
A_{n}\mid a\cdot\underline{1}=\underline{1}\}.
\]

Let $(A_{n})_{n\in\mathbb{N}}\rightarrow(B_{n})_{n\in\mathbb{N}{}}$ be an
inverse systems of injective homomorphisms. From%
\[
0\rightarrow(A_{n})_{n\in\mathbb{N}}\rightarrow(B_{n})_{n\in\mathbb{N}{}%
}\rightarrow(B_{n}/A_{n})_{n\in\mathbb{N}{}}\rightarrow0
\]
we obtain an exact sequence%
\begin{equation}
1\rightarrow\varprojlim A_{n}\rightarrow\varprojlim B_{n}\rightarrow
\varprojlim(B_{n}/A_{n})\rightarrow\varprojlim\nolimits^{1}A_{n}%
\rightarrow\varprojlim\nolimits^{1}B_{n} \label{e0}%
\end{equation}
of groups and pointed sets. Exactness at $\varprojlim(B_{n}/A_{n})$ means that
the fibres of $\varprojlim(B_{n}/A_{n})\rightarrow\varprojlim\nolimits^{1}%
A_{n}$ are the orbits for the natural action of $\varprojlim B_{n}$ on
$\varprojlim(B_{n}/A_{n})$. When each $A_{n}$ is normal in $B_{n}$, so that
$C_{n}=_{\text{df}}B_{n}/A_{n}$ is a group, (\ref{e0}) can be extended to an
exact sequence%
\begin{equation}
1\rightarrow\varprojlim A_{n}\rightarrow\varprojlim B_{n}\rightarrow
\varprojlim C_{n}\rightarrow\varprojlim\nolimits^{1}A_{n}\rightarrow
\varprojlim\nolimits^{1}B_{n}\rightarrow\varprojlim\nolimits^{1}%
C_{n}\rightarrow1. \label{e1}%
\end{equation}
Exactness at $\varprojlim\nolimits^{1}A_{n}$ means that the fibres of
$\varprojlim\nolimits^{1}A_{n}\rightarrow\varprojlim\nolimits^{1}B_{n}$ are
the orbits for the natural action of $\varprojlim C_{n}$ on $\varprojlim
\nolimits^{1}A_{n}$.

Recall that an inverse system $(X_{n})_{n\in\mathbb{N}{}}$ of sets (or groups)
is said to satisfy the condition (ML) if, for each $m$, the decreasing chain
in $X_{m}$ of the images of the $X_{n}$ for $n\geq m$ is eventually constant.

\begin{proposition}
\label{i1}Let $(A_{n},u_{n})_{n\in\mathbb{N}{}}$ be an inverse system of groups.

\begin{enumerate}
\item If $(A_{n},u_{n})$ satisfies (ML), then $\varprojlim\nolimits^{1}%
A_{n}=0$.

\item If the $A_{n}$ are countable and $(A_{n},u_{n})_{n\in\mathbb{N}{}}$
fails (ML), then $\varprojlim\nolimits^{1}A_{n}$ is uncountable.
\end{enumerate}
\end{proposition}

\begin{proof}
(a) The action\footnote{We usually omit the subscript on transition maps.}
\[
(\ldots,a_{n},\ldots)\cdot(\ldots,x_{n},\ldots)=(\ldots,a_{n}\cdot x_{n}%
\cdot(ua_{n+1})^{-1}\text{~},\ldots)
\]
of the group $G_{N+1}=_{\text{df}}%
%TCIMACRO{\tprod \nolimits_{0\leq n\leq N+1}}%
%BeginExpansion
{\textstyle\prod\nolimits_{0\leq n\leq N+1}}
%EndExpansion
A_{n}$ on the set $S_{N}=_{\text{df}}%
%TCIMACRO{\tprod \nolimits_{0\leq n\leq N}}%
%BeginExpansion
{\textstyle\prod\nolimits_{0\leq n\leq N}}
%EndExpansion
A_{n}$ is transitive, and the projection $(a_{n})_{n}\mapsto a_{N+1}$ gives an
isomorphism from the stabilizer of any $x\in$ $S_{N}$ onto $A_{N+1}$. Let
$x,y\in\prod_{n\in\mathbb{N}{}}A_{n}$, and let
\[
P_{N}=\left\{  g\in G_{N+1}\mid gx^{N}=y^{N}\right\}  \text{.}%
\]
where $x^{N}$ and $y^{N}$ are the images of $x$ and $y$ in $S_{N}$. We have to
show that $\varprojlim P_{N}$ is nonempty. The observations in the first
sentence show that $P_{N}$ is nonempty and that $A_{N+1}$ acts simply
transitively on it. It follows that the inverse system $(P_{N})$ satisfies
(ML). Let $Q_{N}=\cap_{i}\mathrm{Im}(P_{N+i}\rightarrow P_{N})$. Then each
$Q_{N}$ is nonempty, and $(Q_{N})_{N\in\mathbb{N}{}}$ is an inverse system
with surjective transition maps. Hence, $\varprojlim Q_{N}$ is (obviously)
nonempty, and any element of it is an element of $\varprojlim P_{N}$.

(b) If $(A_{n})$ fails (ML), then there exists an $m$ such that infinitely
many of the groups%
\[
B_{i}\overset{\text{df}}{=}\text{Im}(A_{m+i}\rightarrow A_{m})
\]
are distinct. As
\[
\varprojlim_{i}{}^{1}A_{m+i}\rightarrow\varprojlim_{i}{}^{1}B_{i}%
\]
is surjective (see (\ref{e1})), it suffices to show that $\varprojlim
\nolimits^{1}B_{i}$ is uncountable. This is accomplished by the next lemma
(applied with $A=A_{m}$).
\end{proof}

\begin{lemma}
\label{i2}Let $\cdots\supset A_{n}\supset A_{n+1}\supset\cdots$ be a sequence
of subgroups of a countable group $A$. If infinitely many of the $A_{n}$ are
distinct, then $\varprojlim\nolimits^{1}A_{n}$ is uncountable.
\end{lemma}

\begin{proof}
From (\ref{e0}) applied to the inverse system $(A_{n}\hookrightarrow
A)_{n\in\mathbb{N}{}}$, we obtain a bijection%
\[
A\backslash(\varprojlim A/A_{n})\rightarrow\varprojlim\nolimits^{1}%
A_{n}\text{.}%
\]
As a map of sets, $A/A_{n+1}\rightarrow A/A_{n}$ is isomorphic to the
projection map
\[
A/A_{n}\times A_{n}/A_{n+1}\rightarrow A/A_{n},
\]
and so $\varprojlim A/A_{n}\approx\prod A_{n}/A_{n+1}$ (as sets), which is uncountable.
\end{proof}

\begin{remark}
\label{i3}The above statements apply to inverse systems indexed by any
directed set $I$ containing an infinite countable cofinal set, because such an
$I$ will also contain a cofinal set isomorphic to $(\mathbb{N}{},\leq)$.
\end{remark}

\begin{note}
\label{i4}The definition of $\varprojlim\nolimits^{1}$ for nonabelian groups
and the sequence (\ref{e1}) can be found in Bousfield and Kan 1972, IX \S 2.
In the commutative case, statement (a) of Proposition \ref{i1} is proved in
Atiyah 1961 and statement (b) in Gray 1966.
\end{note}

\subsection{Torsors}

Let $\mathsf{E}$ be a category with finite fibred products (in particular, a
final object $S$) endowed with a topology in the sense of Grothendieck (see
Bucur and Deleanu 1968, Chapter 2). Thus, $\mathsf{E}$ is a site. By
\textquotedblleft torsor\textquotedblright\ we mean \textquotedblleft right
torsor\textquotedblright.

\begin{plain}
\label{t9}For a sheaf of groups $A$ on $\mathsf{E}$, a right $A$-sheaf $X$,
and a left $A$-sheaf $Y$, $X\wedge^{A}Y$ denotes the contracted product of $X$
and $Y$, i.e., the quotient sheaf of $X\times Y$ by the diagonal action of
$A$, $(x,y)a=(xa,a^{-1}y)$. When $A\rightarrow B$ is a homomorphism of sheaves
of groups, $X\wedge^{A}B$ is the $B$-sheaf obtained from $X$ by extension of
the structure group. In this last case, if $X$ is an $A$-torsor, then
$X\wedge^{A}B$ is a $B$-torsor.
\end{plain}

\begin{plain}
\label{t0}For an $A$-torsor $P$ and a left $A$-sheaf $X$, define%
\[
^{P}X=P\wedge^{A}X.
\]
When $X$ is a sheaf of groups and $A$ acts by group homomorphisms, ${}^{P}X$
is a sheaf of groups. For example, when we let $A$ act on itself by inner
automorphisms, ${}^{P}A$ is the inner form of $A$ defined by $P$. There is a
natural left action of ${}^{P}A$ on $P$, which makes $P$ into a left ${}^{P}%
A$-torsor and induces an isomorphism ${}$%
\begin{equation}
^{P}A\rightarrow\mathcal{A}{}ut_{A}(P). \label{e9}%
\end{equation}
Let $\mathsf{P}(A)$ denote the category of $A$-torsors and $H^{1}(S,A)$ the
set of isomorphism classes of objects in $\mathsf{P}(A)$ (pointed by the class
of the trivial torsor $A_{A}$).
\end{plain}

\begin{plain}
\label{t1}Let $v\colon B\rightarrow C$ be a homomorphism of sheaves of groups
on $\mathsf{E}$, and let $Q$ be a $C$-torsor. Define $\mathsf{P}(B\rightarrow
C;Q)$ to be the category whose objects are the $v$-morphisms of torsors
$P\rightarrow Q$ and whose morphisms $\Hom(P\rightarrow Q,P^{\prime
}\rightarrow Q)$ are the $B$-morphisms $P\rightarrow P^{\prime}$ giving a
commutative triangle. Let $H^{1}(S,B\rightarrow C;Q)$ denote the set of
isomorphism classes in $\mathsf{P}(B\rightarrow C;Q)$. When $Q=C_{C}$, we drop
it from the notation; then $H^{1}(S,B\rightarrow C)$ is pointed by the class
of $B_{B}\rightarrow C_{C}$. The category $\mathsf{P}(B\rightarrow0)$ is
canonically equivalent with the category of $B$-torsors, and so
\[
H^{1}(S,B\rightarrow0)\cong H^{1}(S,B).
\]
Let $A=\Ker(B\rightarrow C)$. Then $A$ is stable under the action of $B$ on
itself by inner automorphisms, and for any object $P\rightarrow Q$ of
$\mathsf{P}(B\rightarrow C;Q)$,
\[
\mathcal{A}{}ut(P\rightarrow Q)\cong P\wedge^{B}A.
\]

\end{plain}

\begin{plain}
\label{t2}Let $v\colon B\rightarrow C$ be a surjective homomorphism with
kernel $A$. To give a $B$-torsor $P$ with $vP$ trivialised by $e\in(vP)(S)$
amounts to giving the $A$-torsor $f^{-1}(e)$: the natural functor
\[
\mathsf{P}(A\rightarrow0)\rightarrow\mathsf{P}(B\rightarrow C)
\]
is an equivalence.
\end{plain}

\begin{plain}
\label{t6}Let $v\colon B\rightarrow C$ be a homomorphism of sheaves of groups
on $\mathsf{E}$. A $B$-torsor $P$ allows us to twist $v$:%
\[
{}^{P}v\colon{}^{P}B\rightarrow{}^{P}C,\quad{}{}^{P}C\overset{\text{df}}%
{=}P\wedge^{B}C.
\]
Here a local section $b$ of $B$ acts on $C$ by $c\mapsto(vb)c(vb)^{-1}$. Let
$Q=vP$. Then ${}^{P}C\cong{}^{Q}C$.
\end{plain}

\begin{plain}
\label{t3}Let $v\colon B\rightarrow C$ be a homomorphism of sheaves of groups
on $\mathsf{E}$. A $B$-torsor $P$ can be regarded as a $({}^{P}B,B)$-bitorsor
(see (\ref{e9})). There is a functor
\begin{equation}
\mathsf{P}({}^{P}B\rightarrow{}{}^{Q}C)\rightarrow\mathsf{P}(B\rightarrow C;Q)
\label{e06}%
\end{equation}
sending $P^{\prime}\rightarrow Q^{\prime}$ to $P^{\prime}\wedge^{{}^{P}%
B}P\rightarrow Q^{\prime}\wedge^{{}^{Q}B}Q$. In particular, the neutral object
of $\mathsf{P}(^{P}B\rightarrow{}^{Q}C)$ is sent to the object $P\rightarrow
Q$ of $\mathsf{P}(B\rightarrow C;Q)$. Let $P^{\text{opp}}$ denote the
$(B,{}^{P}B)$-bitorsor with the same underlying sheaf as $P$ but with local
sections $b$ and $b^{\prime}$ of $B$ and ${}^{P}B$ acting as $(b,b^{\prime
})\cdot p=b^{\prime-1}\cdot p\cdot b^{-1}$. The functor%
\begin{equation}
\mathsf{P}(B\rightarrow C;Q)\rightarrow\mathsf{P}({}^{P}B\rightarrow{}{}^{Q}C)
\label{e07}%
\end{equation}
sending $P^{\prime}\rightarrow Q$ to $(P^{\prime}\rightarrow Q)\wedge
^{B}P^{\text{opp}}$ is a quasi-inverse to the functor in (\ref{e06}).
Therefore, both functors are equivalences of categories.
\end{plain}

\begin{proposition}
\label{t4}Let%
\begin{equation}
1\rightarrow A\rightarrow B\overset{v}{\rightarrow}C\rightarrow0 \label{e03}%
\end{equation}
be an exact sequence of sheaves of groups on $\mathsf{E}$, and let
$P\rightarrow Q$ be a $v$-morphism of torsors. There is a natural bijection%
\[
H^{1}(S,{}^{P}A)\rightarrow H^{1}(S,B\rightarrow C;Q)
\]
sending the neutral element of $H^{1}(S,{}^{P}A)$ to the element
$[P\rightarrow Q]$ of $H^{1}(S,B\rightarrow C;Q)$.
\end{proposition}

\begin{proof}
We can use $P$ to twist the sequence (\ref{e03}):%
\[
1\rightarrow{}^{P}A\rightarrow{}^{P}B\rightarrow{}^{Q}C\rightarrow1,\quad
{}^{P}A=P\wedge^{B}A\text{.}%
\]
Now combine%
\[
H^{1}(S,{}^{P}A)\overset{\ref{t2}}{\rightarrow}H^{1}(S,{}^{P}B\rightarrow
{}^{Q}C)\overset{(\ref{e06})}{\rightarrow}H^{1}(S,B\rightarrow C;Q)\text{.}%
\]

\end{proof}

\begin{remark}
\label{t5}If in the proposition $A$ is commutative, then the action of $B$ on
$A$ factors through an action of $C$ on $A$, and so%
\[
{}^{P}A\overset{\text{df}}{=}P\wedge^{B}A\cong P\wedge^{B}C\wedge^{C}A\cong
Q\wedge^{C}A\overset{\text{df}}{=}{}^{Q}A\text{.}%
\]

\end{remark}

\begin{note}
\label{t7}The basic definitions \ref{t9}--\ref{t0} are from Giraud 1971. The
remaining statements can be found, or are hinted at, in Deligne 1979a,
2.4.3--2.4.4. See also Breen 1990. (The main ideas go back to Dedecker and
Grothendieck in the 1950s.)
\end{note}

\subsection{Cohomology and inverse limits}

We now fix an affine scheme $S$ and let $\mathsf{E}$ be the category of affine
schemes over $S$ endowed with the fpqc topology (that for which the covering
families are the finite surjective families of flat morphisms).

Throughout this subsection, $(G_{n},u_{n})_{n}$ is an inverse system, indexed
by $(\mathbb{N}{},\leq)$, of flat affine group schemes of finite type over $S$
with faithfully flat transition maps, and $G=\varprojlim G_{n}$. Thus, $G$ is
a flat affine group scheme over $S$.

\begin{proposition}
\label{c1}The map $[P]\mapsto([P\wedge^{G}G_{n}])_{n\geq0}$
\begin{equation}
H^{1}(S,G)\rightarrow\varprojlim_{n}H^{1}(S,G_{n}) \label{e4}%
\end{equation}
is surjective. For a $G$-torsor $P$, the fibre of the map containing $[P]$ is
$\varprojlim\nolimits^{1}G_{n}^{\prime}(S)$ where $G_{n}^{\prime}$ is the
inner form ${}P\wedge^{G}G_{n}$ of $G_{n}$.
\end{proposition}

\begin{proof}
A class $c$ in $\varprojlim H^{1}(S,G_{n})$ is represented by an inverse
system%
\[
P_{0}\leftarrow P_{1}\leftarrow\cdots\leftarrow P_{n}\leftarrow\cdots
\]
with $P_{n}$ a $G_{n}$-torsor. The inverse limit of this system is a
$G$-torsor mapping to $c$.

Let $P^{\prime}$ and $P$ be $G$-torsors such that $P_{n}^{\prime}\approx
P_{n}$ for all $n$, and choose isomorphisms $a_{n}\colon P_{n}^{\prime}$
$\rightarrow P_{n}$. Consider
\[
\begin{CD}
P_{n+1}^{\prime} @>a_{n+1}>> P_{n+1}\\
@VV{v^{\prime}}%
V @VV{v}V\\
P_{n}^{\prime} @>a_{n}>> P_{n}.
\end{CD}
\]
There is a unique isomorphism $b_{n}\colon P_{n}^{\prime}\rightarrow P_{n}$
for which the diagram commutes, i.e., such that%
\[
b_{n}\circ v^{\prime}=v\circ a_{n+1}\text{.}%
\]
Let $e_{n}$ be the element of $\Aut(P_{n})$ such that $e_{n}\circ b_{n}=a_{n}%
$; then%
\[
e_{n}\circ v\circ a_{n+1}=a_{n}\circ v^{\prime}\text{.}%
\]

\noindent Replacing $(a_{n})_{n\geq0}$ with $(c_{n}\circ a_{n})_{n\geq0}$
replaces $(e_{n})_{n\geq0}$ with $(c_{n}\cdot e_{n}\cdot u^{\prime}%
c_{n+1}^{-1})_{n\geq0}$ where $u^{\prime}$ is the transition map
$\Aut(P_{n+1})\rightarrow\Aut(P_{n})$. Thus, the class of $(e_{n})_{n\geq0}$
in $\varprojlim\nolimits^{1}\Aut(P_{n})$ is independent of the choice of the
$a_{n}$. Similarly, it depends only on the isomorphism class of $P^{\prime}$.
Therefore, we have a well-defined map from the fibre containing $[P]$ to
$\varprojlim\nolimits^{1}\Aut(P_{n})$, and it is straightforward to check that
it is a bijection. Finally, (\ref{e9}) allows us to replace $\Aut(P_{n})$ with
$G_{n}^{\prime}(S)$.
\end{proof}

\begin{corollary}
\label{c2}When the $G_{n}$ are commutative, there is an exact sequence%
\[
0\rightarrow\varprojlim\nolimits^{1}G_{n}(S)\rightarrow H^{1}(S,G)\rightarrow
\varprojlim H^{1}(S,G_{n})\rightarrow0\text{.}%
\]

\end{corollary}

\begin{proof}
In this case, $G_{n}^{\prime}=G_{n}$.
\end{proof}

\begin{corollary}
\label{c3}For any countable family $(G_{i})_{i\in I}$ of flat affine group
schemes,
\[
H^{1}(S,%
%TCIMACRO{\tprod }%
%BeginExpansion
{\textstyle\prod}
%EndExpansion
G_{i})=%
%TCIMACRO{\tprod }%
%BeginExpansion
{\textstyle\prod}
%EndExpansion
H^{1}(S,G_{i}).
\]

\end{corollary}

\begin{proof}
This is certainly true for finite families. Thus, we may assume that $I$ is
infinite, and equals $\mathbb{N}{}$. Let $A_{n}=\prod_{0\leq i\leq n}G_{i}$.
For any $\prod_{i\geq0}G_{i}$-torsor $P$, the projection maps ${}{}^{P}%
A_{n}(S{})\rightarrow{}{}^{P}A_{n-1}(S)$ admit sections, and so are
surjective. Therefore $\varprojlim\nolimits^{1}{}^{P}A_{n}(S)=0$ (\ref{i1}a),
and it follows that
\[
H^{1}(S{},G)\overset{\ref{c1}}{\cong}\varprojlim_{n}H^{1}(S,A_{n}%
)\cong\varprojlim_{n}%
%TCIMACRO{\tprod \limits_{0\leq i\leq n}}%
%BeginExpansion
{\textstyle\prod\limits_{0\leq i\leq n}}
%EndExpansion
H^{1}(S,G_{i})\cong%
%TCIMACRO{\tprod \limits_{i\geq0}}%
%BeginExpansion
{\textstyle\prod\limits_{i\geq0}}
%EndExpansion
H^{1}(S,G_{i})\text{.}%
\]

\end{proof}

\begin{remark}
\label{c4}Let $S=\Spec(\mathbb{Q}{})$. Although the maps $u_{n}\colon
G_{n}\rightarrow G_{n-1}$ are surjective, typically the maps $G_{n}%
(\mathbb{Q}{})\rightarrow G_{n-1}(\mathbb{Q}{})$ will not be. In fact,
typically, the inverse system $(G_{n}(\mathbb{Q}{}))_{n}$ will not satisfy
(ML) and so $\varprojlim\nolimits^{1}G_{n}(\mathbb{Q}{})$ will be uncountable
(\ref{i1}b). For example, consider a tower of distinct subfields of
$\mathbb{Q}{}^{\text{al}}$,%
\[
\mathbb{Q}{}\subset F_{1}\subset\cdots\subset F_{n-1}\subset F_{n}%
\subset\cdots,\quad\lbrack F_{n}:\mathbb{Q}{}]<\infty\text{.}%
\]
There is an inverse system $(G_{n},u_{n})$ with surjective transition maps for
which $G_{n}$ is the $\mathbb{Q}{}$-torus obtained from $\mathbb{G}_{m/F_{n}}$
by restriction of scalars and $u_{n}$ is the norm map. Then
\[
(G_{n}(\mathbb{Q}{}),u_{n}(\mathbb{Q}{}))_{n\in\mathbb{N}{}}=(F_{n}^{\times
},\Nm_{F_{n}/F_{n-1}})_{n\in\mathbb{N}{}}\text{,}%
\]
which fails (ML),\footnote{To see this, use that, for a finite extension $E/F$
of number fields and a finite prime $v$ of $F$, $\ord_{v}(\Nm E^{\times})$ is
the ideal in $\mathbb{Z}{}$ generated by the residue class degrees of the
primes of $E$ lying over $v$.} and so $\varprojlim\nolimits^{1}G_{n}%
(\mathbb{Q}{})$ is uncountable.
\end{remark}

\subsection{Comparison with Galois cohomology}

We now let $S$ be the spectrum of a field $k$, and we let $H^{1}(k,-)$ denote
$H^{1}(S,-)$. Choose a separable closure $k^{\text{sep}}$ of $k$, and let
$\Gamma=\Gal  (k^{\text{sep}}/k)$.

\begin{proposition}
\label{c5}For any smooth algebraic group $N$ over $k$, there is a canonical
isomorphism%
\[
H^{1}(k,N)\rightarrow H^{1}(\Gamma,N(k^{\text{sep}})).
\]

\end{proposition}

\begin{proof}
An $N$-torsor $P$ is represented by an algebraic variety over $k$, and hence
acquires a point $p$ over some subfield of $k^{\text{sep}}$ of finite degree
over $k$. The formula%
\begin{equation}
\tau p=p\cdot a_{\tau} \label{e5}%
\end{equation}
defines a continuous crossed homomorphism $a_{\tau}\colon\Gamma\rightarrow
N(k^{\text{sep}})$ whose cohomology class is independent of the choice of $p$
and depends only on the isomorphism class of $P$. Thus, we have a well-defined
map $H^{1}(k,N)\rightarrow H^{1}(\Gamma,N(k^{\text{sep}}))$, and it follows
from descent theory that this is an isomorphism.
\end{proof}

\begin{plain}
\label{c11}Let $N$ be a smooth algebraic group over $k$, and let
$f\colon\Gamma\rightarrow\mathcal{A}{}ut(N)(k^{\text{sep}})$ be a continuous
crossed homomorphism (discrete topology on $\mathcal{A}{}ut(N)(k^{\text{sep}%
})$). The \textquotedblleft twist\textquotedblright\ of $N$ by $f$ is a smooth
algebraic group ${}_{f}N$ over $k$ such that ${}_{f}N(k^{\text{sep}%
})=N(k^{\text{sep}})$ but with $\tau\in\Gamma$ acting according to the rule%
\[
\tau\ast x=f(\tau)\cdot\tau x\text{.}%
\]
When we let $N$ act on itself by inner automorphisms, a crossed homomorphism
$f\colon\Gamma\rightarrow N(k^{\text{sep}})$ defines a twist $_{f}N$ of $N$
with $\tau\in\Gamma$ acting on $_{f}N(k^{\text{sep}})$ by%
\[
\tau\ast x=f(\tau)\cdot\tau x\cdot f(\tau)^{-1}\text{.}%
\]

\end{plain}

\begin{plain}
\label{c12}Let $G=\varprojlim(G_{n},u_{n})$ be as in the preceding subsection
but with $G_{n}$ now a smooth algebraic group over $k$, and define
$H_{\text{cts}}^{1}(\Gamma,G)$ be the cohomology set computed using crossed
homomorphisms $\Gamma\rightarrow G(k^{\text{sep}})$ that are continuous for
the profinite topology on $\Gamma$ and the inverse limit topology on
$G(k^{\text{sep}})=\varprojlim G_{n}(k^{\text{sep}})$ (discrete topology on
$G_{n}(k^{\text{sep}})$). Thus, giving a continuous crossed homomorphism
$f\colon\Gamma\rightarrow G(k^{\text{sep}})$ amounts to giving a compatible
family of continuous crossed homomorphisms $f_{n}\colon\Gamma\rightarrow
G_{n}(k^{\text{sep}})$.
\end{plain}

\begin{proposition}
\label{c6}The map%
\[
H_{\text{cts}}^{1}(\Gamma,G)\rightarrow\varprojlim H^{1}(\Gamma,G_{n})
\]
sending $[f]$ to $([f_{n}])_{n\geq0}$ is surjective. The fibre of the map
containing $[f]$ equals $\varprojlim\nolimits^{1}{}G_{n}^{\prime}(k)$ where
$G_{n}^{\prime}={}_{f}G_{n}$.
\end{proposition}

\begin{proof}
Each class $c$ in $\varprojlim H^{1}(S,G_{n})$ is represented by a family
$(f_{n})_{n\geq0}$ of crossed homomorphisms, which can be chosen so that
$f_{n-1}=u_{n}\circ f_{n}$. The $f_{n}$ define a continuous crossed
homomorphism $f\colon\Gamma\rightarrow G(k^{\text{sep}})$ mapping to $c$.

Let $f^{\prime}$ and $f$ be continuous crossed homomorphisms such that
$f_{n}^{\prime}\sim f_{n}$ for all $n$, and choose $a_{n}\in G_{n}%
(k^{\text{sep}})$ in such a way that%
\begin{equation}
f^{\prime}(\tau)_{n}=a_{n}^{-1}\cdot f(\tau)_{n}\cdot\tau a_{n}\text{.}
\label{e7}%
\end{equation}
Define $e_{n}\in G_{n}(k^{\text{sep}})$ by the equation%
\[
e_{n}\cdot ua_{n+1}=a_{n}\text{.}%
\]

On applying $u$ to the equation (\ref{e7})$_{n+1}$, we obtain the equation%
\[
f^{\prime}(\tau)_{n}=(ua_{n+1})^{-1}\cdot f(\tau)_{n}\cdot\tau(ua_{n+1}%
)\text{,}%
\]
or,%
\[
f^{\prime}(\tau)_{n}=a_{n}^{-1}\cdot e_{n}\cdot f(\tau)_{n}\cdot\tau
e_{n}^{-1}\cdot\tau a_{n}.
\]
On comparing this with (\ref{e7}), we find that%
\[
e_{n}=f_{n}(\tau)\cdot\tau e_{n}\cdot f_{n}(\tau)^{-1},
\]
ie., that%
\[
e_{n}\in(_{f}G_{n})(k^{\text{sep}})^{\Gamma}=(_{f}G_{n})(k)\text{.}%
\]

The element $a_{n}$ can be replaced by $c_{n}\cdot a_{n}$, where $c_{n}$ is
any element of $(_{f}G_{n})(k)$. When this is done for each $n$,
$(e_{n})_{n\geq0}$ is replaced by $(c_{n}\cdot e_{n}\cdot(uc_{n+1}%
)^{-1})_{n\geq0}$.Thus, the class of $(e_{n})$ in $\varprojlim\nolimits^{1}%
(_{f}G_{n})(k)$ is independent of the choice of the $a_{n}$. Similarly, it
depends only on the cohomology class of $f^{\prime}$. Therefore, we have a
well-defined map from the fibre containing $[f]$ to $\varprojlim
\nolimits^{1}(_{f}G_{n})(k)$, and it is straightforward to check that this is
a bijection.
\end{proof}

\begin{proposition}
\label{c7}There is a canonical isomorphism of pointed sets%
\[
H^{1}(k,G)\rightarrow H_{\text{cts}}^{1}(\Gamma,G).
\]

\end{proposition}

\begin{proof}
Let $P$ be a $G$-torsor, and let $P_{n}=u_{n}P$. Then $P(k^{\text{sep}%
})=\varprojlim P_{n}(k^{\text{sep}})$, which, because the maps $P_{n+1}%
(k^{\text{sep}})\rightarrow P_{n}(k^{\text{sep}})$ are surjective, is
nonempty. Choose a $p\in P(k^{\text{sep}})$. Then the formula%
\[
\tau p=p\cdot f(\tau)
\]
defines a continuous crossed homomorphism $f\colon\Gamma\rightarrow
G(k^{\text{sep}})$ whose cohomology class is independent of the choice of $p$
and of the choice of $P$ in its isomorphism class. Therefore, we have a
well-defined map $H^{1}(k,G)\rightarrow H_{\text{cts}}^{1}(\Gamma,G)$. Since
this map is compatible with the maps in Propositions \ref{c1} and \ref{c6},
they, together with (\ref{c5}) show that it is a bijection.
\end{proof}

\begin{remark}
\label{c9}Let $S$ be the spectrum of a field. The following conditions on an
affine group scheme $G$ over $S$ are equivalent:

\begin{enumerate}
\item the set of closed normal subgroup schemes $H\subset G$ such that $G/H$
is of finite type over $S$ is infinite and countable;

\item $G=\varprojlim_{n}(G_{n},u_{n})$ where $(G_{n},u_{n})_{n}$ is an inverse
system, indexed by $(\mathbb{N}{},\leq)$, of affine algebraic groups over $S$
with surjective transition maps.
\end{enumerate}

\noindent An affine group scheme $G$ satisfying these conditions will be said
to be \emph{separable}. When $G$ is separable, any inverse system satisfying
(b) is cofinal in the inverse system of all algebraic quotients of $G$.
Therefore, the former can be replaced by the latter, which makes Proposition
\ref{c6} more canonical.
\end{remark}

\begin{note}
\label{c10}Proposition \ref{c5} is a standard result. In the commutative case,
Proposition \ref{c6} is proved in Tate 1976.
\end{note}

\subsection{Application to periods}

Let $\Mot(\mathbb{Q}{})$ denote the category of motives based on all smooth
projective varieties over $\mathbb{Q}{}$. Let $\CM(\mathbb{Q}{})$ be the
Tannakian subcategory of $\Mot(\mathbb{Q}{})$ generated by the
zero-dimensional varieties over $\mathbb{Q}{}$ and the abelian varieties of
CM-type (see \S 2). Let $G^{\Mot}=\mathcal{A}{}ut^{\otimes}(\omega_{\text{B}%
})$ and $P^{\Mot}=\mathcal{I}{}som^{\otimes}(\omega_{\text{B}},\omega
_{\text{dR}})$, and define $G^{\CM}$ and $P^{\CM}$ similarly. From the
inclusion $\CM(\mathbb{Q})\subset\Mot(\mathbb{Q}{})$, we obtain a faithfully
flat homomorphism $G^{\Mot}\rightarrow G^{\CM}$.

\begin{theorem}
\label{c8}If the kernel of $G^{\Mot}\rightarrow G^{\CM}$ is an inverse limit
of simply connected semisimple groups, then the isomorphism class of
$P^{\Mot}\rightarrow P^{\CM}$ in $\mathsf{P}(G^{\Mot}\rightarrow
G^{\CM};P^{\CM})$ is uniquely determined by its class over $\mathbb{R}{}$.
\end{theorem}

\begin{proof}
Let $G=\Ker(G^{\Mot}\rightarrow G^{\CM})$. The condition on $G$ implies that
it is, in fact, a product of semisimple groups each of which is simply
connected. Moreover, the product is countable because $\Mot(\mathbb{Q}{})$ is
generated as a Tannakian category by a countable set of varieties. According
to Proposition \ref{t4}, the isomorphism classes in $\mathsf{P}(G^{\Mot}%
\rightarrow G^{\CM};P^{\CM})$ are classified by $H^{1}(\mathbb{Q}{}%
,{}G^{\prime}$) where $G^{\prime}=P^{\Mot}\wedge^{G^{\Mot}}G$. As $G^{\prime}$
is a form of $G$, it also is a countable product of simply connected
semisimple groups, and so the proposition follows from Corollary \ref{c3} and
the theorem of Kneser, Harder, and Chernousov (see \ref{g02} below).
\end{proof}

\begin{remark}
\label{a1}(a) The group $G^{\Mot}$ is proreductive because $\Mot(\mathbb{Q}%
{})$ is semisimple. If Deligne's hope\ (1979b, 0.10) that every Hodge class is
an absolute Hodge class is true, then $(G^{\Mot})^{\circ}$ is the group
attached to the category of motives over $\mathbb{Q}{}^{\text{al}}$ (see
Deligne and Milne 1982, 6.22, 6.23, where the hypothesis was inadvertently
omitted); moreover, the group $G$ in the above proof is the kernel of the
canonical homomorphism from $(G^{\Mot})^{\circ}$ to the Serre group (ibid.
p220), and it is the derived group of $(G^{\Mot})^{\circ}$; it is therefore an
inverse limit of semisimple groups.

(b) It is generally hoped that the derived group of $(G^{\Mot})^{\circ}$ is
simply connected --- see the question in Serre 1994, 8.1.
\end{remark}

\begin{remark}
\label{a2}Let $G=\Ker(G^{\Mot}\rightarrow G^{\CM})$, and let $H\subset G$ be
the intersection of the kernels of the homomorphisms from $G$ onto simply
connected semisimple algebraic groups over $\mathbb{Q}{}$. Let $\Mot^{H}%
(\mathbb{Q}{})$ be the subcategory of $\Mot(\mathbb{Q}{})$ of objects on which
$H$ acts trivially. Then, with the obvious notations, the isomorphism class of
$P^{\Mot^{H}}\rightarrow P^{\CM}$ in $\mathsf{P}(G^{\Mot^{H}}\rightarrow
G^{\CM};P^{\CM})$ is uniquely determined by its class over $\mathbb{R}{}$. The
preceding remark indicates that it is reasonable to hope that $H=0$.
\end{remark}

\newpage

\section{Periods of abelian varieties with complex multiplication}

\subsection{Periods of zero-dimensional varieties}

The category of motives based on the zero-dimensional varieties over
$\mathbb{Q}{}$ is denoted $\Art(\mathbb{Q}{})$ and is called the category of
Artin motives over $\mathbb{Q}$. The Betti fibre functor $\omega_{B}$ defines
an equivalence of $\Art(\mathbb{Q}{})$ with the category of continuous
representations of $\Gamma$ on finite-dimensional $\mathbb{Q}{}$-vector spaces
(Deligne and Milne 1982, 6.17) from which it follows that $\mathcal{A}%
{}ut^{\otimes}(\omega_{B})$ is the constant profinite group scheme $\varGamma$
with $\varGamma(\mathbb{Q}{})=\varGamma(\mathbb{Q}{}^{\text{al}})=\Gamma$. The
de Rham fibre functor is $hX\mapsto\Gamma(X,\mathcal{O}{}_{X})$, from which
the next statement follows easily.

\begin{theorem}
\label{p1}Let $P^{\Art}$ be the period torsor for $\Art(\mathbb{Q}{})$. Then
$P^{\Art}\cong\Spec\mathbb{Q}{}^{\text{al}}$ with its natural action of
$\varGamma$, and $p^{\Art}$ is the obvious $\mathbb{Q}{}^{\text{al}}$-point of
$\Spec\mathbb{Q}{}^{\text{al}}$. Thus, the period point $p^{\Art}$ has
coordinates in $\mathbb{Q}{}^{\text{al}}$, and the cocycle corresponding to
the pair $(P^{\Art},p^{\Art})$ is the crossed homomorphism $\Gamma
\rightarrow\varGamma(\mathbb{Q}{}^{\text{al}})$, $\tau\mapsto\tau$.
\end{theorem}

Note that the theorem determines the pair $(P^{\Art},p^{\Art})$ uniquely up to
a unique isomorphism.

Throughout this section, we use $f$ to denote the (crossed) homomorphism in
the theorem.

\subsection{Notations for tori}

For a finite \'{e}tale $\mathbb{Q}{}$-algebra $A$, let
\begin{align*}
\varSigma_{A/\mathbb{Q}{}}  &  =\Hom_{\mathbb{Q}{}\text{-alg}}(A,\mathbb{Q}%
{}^{\text{al}}),\\
\Gamma_{A/\mathbb{Q}{}}  &  =\Aut_{\mathbb{Q}{}\text{-alg}}(A),
\end{align*}
and let $(\mathbb{G}_{m})_{A/\mathbb{Q}{}}$ be the torus over $\mathbb{Q}{}$
obtained from $\mathbb{G}_{m/A}$ by restriction of scalars. Thus%
\[
(\mathbb{G}_{m})_{A/\mathbb{Q}{}}(R)=(A\otimes R)^{\times}%
\]
for all $\mathbb{Q}{}$-algebras $R$. For an infinite field extension $K/k$,%
\[
(\mathbb{G}_{m})_{K/k}=\varprojlim(\mathbb{G}_{m})_{K^{\prime}/\mathbb{Q}{}%
},\quad K^{\prime}\subset K,\quad\lbrack K^{\prime}:\mathbb{Q}{}%
]<\infty\text{.}%
\]

\begin{plain}
\label{tori1}There is an equivalence $T\mapsto X^{\ast}(T)=_{\text{df}%
}\Hom(T_{/\mathbb{Q}{}^{\text{al}}},\mathbb{G}_{m})$ from the category of tori
over $\mathbb{Q}{}$ to the category of finitely generated free $\mathbb{Z}{}%
$-modules endowed with a continuous left action of $\Gamma$.
\end{plain}

\begin{plain}
\label{tori2}There is an equivalence $A\mapsto\varSigma_{A/\mathbb{Q}{}}$ from
the category of finite \'{e}tale $\mathbb{Q}{}$-algebras to the category of
finite sets endowed with a continuous left action of $\Gamma$. A quasi-inverse
is provided by%
\[
\varSigma\mapsto A(\varSigma)\overset{\text{df}}{=}\Hom_{\Gamma}%
(\varSigma,\mathbb{Q}{}^{\text{al}}).
\]
If $\varSigma\leftrightarrow A$, then the decomposition $\varSigma=\amalg
\varSigma_{i}$ of $\varSigma$ into orbits corresponds to the decomposition
$A=\prod A(\varSigma_{i})$ of $A$ into a product of fields. If $\Gamma$ acts
transitively on $\varSigma$, then the choice of an $e\in\varSigma$ determines
isomorphisms $\Gamma/\Gamma_{e}\rightarrow\varSigma$ and
$A(\varSigma)\rightarrow(\mathbb{Q}{}^{\text{al}})^{\Gamma_{e}}$. Here
$\Gamma_{e}=\{\tau\in\Gamma\mid\tau e=e\}$.
\end{plain}

\begin{plain}
\label{tori3}On combining these equivalences, we see that there is a fully
faithful functor $\varSigma\mapsto T^{\varSigma}$ from the category of finite
sets endowed with a continuous left $\Gamma$-action to the category of tori,
for which $T^{\Sigma}=(\mathbb{G}_{m})_{A(\varSigma)/\mathbb{Q}{}}$, $X^{\ast
}(T^{\Sigma})=\mathbb{Z}{}[\varSigma]$ (free $\mathbb{Z}{}$-module on
$\varSigma$ with the natural left action of $\Gamma$), and $T^{\varSigma_{1}%
\sqcup\varSigma_{2}}=T^{\varSigma_{1}}\times T^{\varSigma_{2}}$.
\end{plain}

\subsection{Abelian varieties with complex multiplication}

For an abelian variety $A$ defined over a subfield $k$ of $\mathbb{C}{}$, the
\emph{Mumford-Tate group} of $A$ is
\[
\MT(A)=\mathcal{A}{}ut^{\otimes}(\omega_{\text{B}}|\langle A_{/\mathbb{C}{}%
}\rangle^{\otimes})\text{.}%
\]
where $\langle A_{/\mathbb{C}{}}\rangle^{\otimes}$ be the category of motives
based on $A_{/\mathbb{C}{}}$ and the projective spaces. There is a canonical
cocharacter $\mu^{A}$ of $\MT(A)$ that splits the Hodge filtration on
$H_{\text{B}}^{1}(A)$. When $\MT(A)$ is commutative (hence a torus), we say
that $A$ is of \emph{CM-type}\footnote{Other authors, and this author at other
times, say that $A$ is potentially of CM-type.}, and we define the
\emph{reflex field }of $A$ to be the field of definition of $\mu^{A}$. It is a
CM-subfield of $\mathbb{C}{}$.

For subfields $k$, $K$ of $\mathbb{C}{}$ with $K$ a CM-field, we define
$\CM^{K}(k)$ to be the category of motives based on

\begin{itemize}
\item the abelian varieties of CM-type over $k$ with reflex field contained in
$K$,

\item the projective spaces, and

\item the zero-dimensional varieties.
\end{itemize}

\noindent When $K=\mathbb{Q}{}^{\text{cm}}$, the composite of all CM-subfields
of $\mathbb{Q}{}^{\text{al}}$, we omit it from the notation.

\subsection{The category $\CM(\mathbb{\mathbb{C}{})}$}

\begin{plain}
\label{cm1}Let $K$ be a CM-subfield of $\mathbb{C}$. With the notation of
(\ref{tori3}), $(\mathbb{G}_{m})_{K/\mathbb{Q}}=T^{\varSigma_{K/\mathbb{Q}{}}%
}$ and $X^{\ast}((\mathbb{G}_{m})_{K/\mathbb{Q}})=\mathbb{Z}%
[\varSigma_{K/\mathbb{Q}{}}]$. The group $S^{K}=\mathcal{A}{}ut^{\otimes
}(\omega_{B}|\CM^{K}(\mathbb{C}))$ is called the \emph{Serre group} for $K$.
When $K$ has finite degree over $\mathbb{Q}{}$, $S^{K}$ is the quotient of
$(\mathbb{G}_{m})_{K/\mathbb{Q}}$ such that%
\[
X^{\ast}(S^{K})=\{n\in\mathbb{Z}{}[\varSigma_{K/\mathbb{Q}{}}]\mid n+\iota
n=\text{constant}\}\text{.}%
\]
Thus, there is an exact sequence%
\begin{equation}
0\rightarrow(\mathbb{G}_{m})_{F/\mathbb{Q}{}}\rightarrow(\mathbb{G}%
_{m})_{K/\mathbb{Q}{}}\times\mathbb{G}_{m}\rightarrow S^{K}\rightarrow0
\label{e04}%
\end{equation}
where $F$ is the largest totally real subfield of $K$. We have
\[
S=\varprojlim_{K\in\mathcal{K}{}}S^{K}%
\]
where $\mathcal{K}{}$ is the set of all CM-subfields of $\mathbb{C}$ finite
and Galois over $\mathbb{Q}{}$.
\end{plain}

\begin{plain}
\label{cm2}For an abelian variety $A$ over $\mathbb{C}{}$, there is a
canonical surjection $S\rightarrow\MT(A)$, which factors through $S^{K}$ if
and only if $A$ has reflex field contained in $K$.
\end{plain}

\begin{note}
For more on Mumford-Tate groups, see Deligne 1982a, \S 3.
\end{note}

\subsection{The category $\CM(\mathbb{Q}^{\text{al}})$}

\begin{plain}
\label{cm3}The functor the category of abelian varieties of CM-type up to
isogeny over $\mathbb{Q}{}^{\text{al}}$ to the similar category over
$\mathbb{C}{}$ is an equivalence. For any CM-subfield $K$ of $\mathbb{C}{}$,
the functor \textquotedblleft extension of scalars\textquotedblright%
\ $\CM^{K}(\mathbb{Q}^{\text{al}})\rightarrow\CM^{K}(\mathbb{C})$ is an
equivalence of tensor categories.
\end{plain}

\subsection{The period torsor}

\begin{plain}
\label{p0}Let $K$ be a CM-subfield of $\mathbb{C}{}$, Galois over
$\mathbb{Q}{}$. Then $\Gamma$ acts on the terms of the sequence (\ref{e04}) in
such a way that the action on the sequence of $\mathbb{Q}{}$-points%
\[
1\rightarrow F^{\times}\rightarrow K^{\times}\times\mathbb{Q}{}^{\times
}\rightarrow S^{K}(\mathbb{Q}{})\rightarrow1
\]
is the obvious one. Let
\[
G_{\text{B}}^{\CM,K}=\mathcal{A}{}ut^{\otimes}(\omega_{\text{B}}%
|\CM^{K}(\mathbb{Q})).
\]
The tensor functors%
\begin{equation}
\CM^{K}(\mathbb{Q}^{\text{al}})\leftarrow\CM^{K}(\mathbb{Q})\leftarrow
\Art(\mathbb{Q}{}) \label{e13}%
\end{equation}
define homomorphisms%
\begin{equation}
1\rightarrow S^{K}\rightarrow G_{\text{B}}^{\CM,K}\overset{v}{\rightarrow
}\varGamma\rightarrow1\text{.} \label{e12}%
\end{equation}
This sequence is exact and the action of $\varGamma$ on $S^{K}$ it defines is
that described above (Deligne 1982b, Lemme 1, Lemme 2).
\end{plain}

\begin{plain}
Let $P^{\CM,K}$ be the period torsor for $\CM^{K}(\mathbb{Q})$. The second
functor in (\ref{e13}) defines a $v$-morphism $P^{\CM,K}\rightarrow
P^{\text{Art}}$, and since $P^{\text{Art}}$ is known (Theorem \ref{p1}), in
order to determine $P^{\CM,K}$ it suffices to determine $P^{\CM,K}\rightarrow
P^{\text{Art}}$ as an object of $\mathsf{P}(G_{\text{B}}^{\CM,K}%
\rightarrow\varGamma;P^{\text{Art}})$. The first step is to correctly identify
the cohomology group classifying the isomorphism classes of objects in this category.
\end{plain}

\begin{proposition}
\label{p00}Let $f\colon\Gamma\rightarrow\varGamma(\mathbb{Q}{})$ be as in
Theorem \ref{p1}, and let $_{f}S^{K}$ be the twist of $S^{K}$ by $f$. There is
a natural one-to-one correspondence between the set of isomorphism classes in
$\mathsf{P}(G_{\text{B}}^{\CM,K}\rightarrow\varGamma;P^{\Art})$ and
$H^{1}(\mathbb{Q},{}_{f}S^{K})$ (fpqc cohomology).
\end{proposition}

\begin{proof}
Apply Proposition \ref{t4}, Remark \ref{t5}, and Theorem \ref{p1}.
\end{proof}

In particular, when $K=\mathbb{Q}^{\text{cm}}$, the isomorphism classes of
objects $P\rightarrow P^{\Art}$ are classified by the group $H^{1}%
(\mathbb{Q}{},{}_{f}S)$. According to Corollary \ref{c2}, there is an exact
sequence%
\[
0\rightarrow\varprojlim\nolimits_{K\in\mathcal{K}{}}^{1}{}_{f}S^{K}%
(\mathbb{Q}{})\rightarrow H^{1}(\mathbb{Q},{}_{f}S)\rightarrow\varprojlim
\nolimits_{K\in\mathcal{K}{}}H^{1}(\mathbb{Q},{}_{f}S^{K})\rightarrow0\text{.}%
\]
We shall show:

\begin{proposition}
\label{p01}The group $\varprojlim\nolimits_{K\in\mathcal{K}{}}H^{1}%
(\mathbb{Q},{}_{f}S^{K})=0$, but $\varprojlim\nolimits_{K\in\mathcal{K}{}}%
^{1}{}_{f}S^{K}(\mathbb{Q}{})$ is uncountable. Therefore,
\[
H^{1}(\mathbb{Q},{}_{f}S)\cong\varprojlim\nolimits_{K\in\mathcal{K}{}}^{1}%
{}_{f}S^{K}(\mathbb{Q}{})
\]
and is uncountable.
\end{proposition}

\subsubsection{Twisting $(\mathbb{G}_{m})_{L/\mathbb{Q}{}}$}

\begin{plain}
\label{cm4}Let $L$ be a finite extension of $\mathbb{Q}{}$, and let
$T=(\mathbb{G}_{m})_{L/\mathbb{Q}{}}$.

\begin{enumerate}
\item The (left) action of $\tau\in\Gamma$ on $\varSigma_{L/\mathbb{Q}{}}$,
$\tau\sigma=\tau\circ\sigma$, corresponds to the natural (left) action of
$\Gamma$ on $X^{\ast}(T)=\mathbb{Z}{}[\varSigma_{L/\mathbb{Q}{}}]$ (see
\ref{tori3}).

\item Assume $L\subset\mathbb{Q}{}^{\text{al}}$ and is Galois over
$\mathbb{Q}{}$. Identify $\varSigma_{L/\mathbb{Q}{}}$ with $\Gamma
_{L/\mathbb{Q}{}}$. There is then a natural (left) action of $\Gamma$ on $T$
(as a torus over $\mathbb{Q}{}$), which defines a (left) action of $\Gamma$ on
$X^{\ast}(T)=\mathbb{Z}{}[\Gamma_{L/\mathbb{Q}{}}]$. The former gives the
obvious action of $\Gamma$ on $T(\mathbb{Q})=L^{\times}$, and the latter
corresponds to the (left) action of $\tau\in\Gamma$ on $\Gamma_{L/\mathbb{Q}%
{}}$, $\tau\sigma=\sigma\circ\tau^{-1}$.
\end{enumerate}
\end{plain}

\begin{lemma}
\label{p03}For a subfield $L$ of $\mathbb{Q}{}^{\text{al}}$, finite and Galois
over $\mathbb{Q}{}$, let $_{f}\Gamma_{L/\mathbb{Q}{}}$ denote the twist of
$\Gamma_{L/\mathbb{Q}{}}$ by $f$ (so that $\tau\in\Gamma$ acts by $\tau
\sigma=\tau\circ\sigma\circ\tau^{-1}$ (Serre 1964, 5.3)), and let
\[
B(L)=A(_{f}\Gamma_{L/\mathbb{Q}{}})\text{.}%
\]
Then%
\[
{}_{f}(\mathbb{G}_{m})_{L/\mathbb{Q}{}}\cong(\mathbb{G}_{m})_{B(L)/\mathbb{Q}%
{}}.
\]

\end{lemma}

\begin{proof}
Clearly, $_{f}(\mathbb{G}_{m})_{L/\mathbb{Q}{}}\cong T^{{}_{f}\Gamma
_{L/\mathbb{Q}{}}}$, which equals $(\mathbb{G}_{m})_{A(_{f}\Gamma
_{L/\mathbb{Q}{}})/\mathbb{Q}{}}$.
\end{proof}

Note that the orbits of $\Gamma$ acting on $_{f}\Gamma_{L/\mathbb{Q}{}}$ are
the conjugacy classes $C$ in $\Gamma_{L/\mathbb{Q}{}}$, and so (\ref{tori2},
\ref{tori3})%
\[
B(L)=%
%TCIMACRO{\tprod _{C}}%
%BeginExpansion
{\textstyle\prod_{C}}
%EndExpansion
A(C),\quad(\mathbb{G}_{m})_{B(L)/\mathbb{Q}{}}=%
%TCIMACRO{\tprod _{C}}%
%BeginExpansion
{\textstyle\prod_{C}}
%EndExpansion
T^{C}\text{.}%
\]
For any $\sigma\in C$, $A(C)\cong L^{Z(\sigma)}$ where $Z(\sigma)$ is the
centralizer of $\sigma$ in $\Gamma_{L/\mathbb{Q}{}}$.

\subsubsection{The tori $_{f}S^{K}$ and $_{f}\bar{S}^{K}$}

Let $K\in\mathcal{K}$. Let $w\colon\mathbb{G}_{m}\rightarrow S^{K}$ be the
weight homomorphism, and let $\bar{S}^{K}=S^{K}/w(\mathbb{G}_{m})$. Then%
\[
X^{\ast}(\bar{S}^{K})=\{n\in\mathbb{Z}{}[\varSigma_{K}]\mid n+\iota
n=1\}\text{, }%
\]
and there is an exact sequence%
\begin{equation}
1\rightarrow(\mathbb{G}_{m})_{F/\mathbb{Q}{}}\rightarrow(\mathbb{G}%
_{m})_{K/\mathbb{Q}{}}\rightarrow\bar{S}^{K}\rightarrow1\text{.} \label{e01}%
\end{equation}
The group $\Gamma$ acts on the exact sequence (\ref{e01}), and so we can twist
the sequence by $f$ to obtain an exact sequence%
\begin{equation}
1\rightarrow{}_{f}(\mathbb{G}_{m})_{F/\mathbb{Q}{}}\rightarrow{}%
_{f}(\mathbb{G}_{m})_{K/\mathbb{Q}{}}\rightarrow{}_{f}\bar{S}^{K}%
\rightarrow1\text{.} \label{e02}%
\end{equation}

\begin{lemma}
\label{p4}If $K$ contains a quadratic imaginary field $k$, then%
\[
{}_{f}\bar{S}^{K}\approx({}\mathbb{G}_{m})_{B(F)/\mathbb{Q}{}}%
\]
where $F$ is the largest totally real subfield of $K$.
\end{lemma}

\begin{proof}
As%
\[
\Gamma_{K/\mathbb{Q}{}}=\Gamma_{F/\mathbb{Q}{}}\times\Gamma_{k/\mathbb{Q}{}%
},\quad\Gamma_{k/\mathbb{Q}{}}=\{1,\iota\},
\]
for each conjugacy class $C$ in $\Gamma_{F/\mathbb{Q}{}}$, there are exactly
two conjugacy classes in $\Gamma_{K/\mathbb{Q}{}}$ mapping to it, namely,%
\[
C_{1}=\{(\tau,1)\mid\tau\in C\},\quad C_{\iota}=\{(\tau,\iota)\mid\tau\in
C\}\text{.}%
\]
Therefore, (\ref{e02}) can be written%
\[
1\rightarrow%
%TCIMACRO{\tprod \nolimits_{C}}%
%BeginExpansion
{\textstyle\prod\nolimits_{C}}
%EndExpansion
T^{C}\rightarrow%
%TCIMACRO{\tprod \nolimits_{C}}%
%BeginExpansion
{\textstyle\prod\nolimits_{C}}
%EndExpansion
(T^{C_{1}}\times T^{C_{\iota}})\rightarrow{}{}_{f}\bar{S}^{K}\rightarrow1
\]
(product over the conjugacy classes in $\Gamma_{F/\mathbb{Q}{}}$). Since%
\[
C\cong C_{1}\cong C_{\iota}\quad\text{(as }\Gamma\text{-sets),}%
\]
each of the maps
\[%
%TCIMACRO{\tprod \nolimits_{C}}%
%BeginExpansion
{\textstyle\prod\nolimits_{C}}
%EndExpansion
T^{C}\rightarrow%
%TCIMACRO{\tprod \nolimits_{C}}%
%BeginExpansion
{\textstyle\prod\nolimits_{C}}
%EndExpansion
T^{C_{1}},\quad%
%TCIMACRO{\tprod \nolimits_{C}}%
%BeginExpansion
{\textstyle\prod\nolimits_{C}}
%EndExpansion
T^{C_{1}}\rightarrow{}_{f}\bar{S}^{K}%
\]
is an isomorphism.
\end{proof}

\medskip

\noindent\textbf{The group $\varprojlim_{K\in\mathcal{K}{}}H^{1}(\mathbb{Q}%
,{}_{f}S^{K})$.}

\begin{proposition}
[Wintenberger 1990, 1.3]\label{p5}If $K$ contains a quadratic imaginary number
field, then $H^{1}(\mathbb{Q}{},{}_{f}S^{K})=0$.
\end{proposition}

\begin{proof}
From the exact sequence%
\begin{equation}
1\rightarrow\mathbb{G}_{m}\rightarrow{}_{f}S^{K}\rightarrow{}_{f}\bar{S}%
^{K}\rightarrow1 \label{e05}%
\end{equation}
we see that it suffices to show that $H^{1}(\mathbb{Q}{},{}_{f}\bar{S}^{K}%
)=0$, but this follows from Lemma \ref{p4} and Hilbert's Theorem 90.
\end{proof}

\begin{corollary}
\label{p6}The group $\varprojlim_{K\in\mathcal{K}{}}H^{1}(\mathbb{Q}{},{}%
_{f}\bar{S}^{K})=0$.
\end{corollary}

\begin{proof}
The CM-fields satisfying the hypothesis of the proposition are cofinal in
$\mathcal{K}{}$.
\end{proof}

\begin{remark}
\label{p7}Wintenberger (1990, p3) shows by example that the proposition is
false without the hypothesis on $K$.
\end{remark}

\medskip

\noindent\textbf{The group $\varprojlim\nolimits_{K\in\mathcal{K}{}}^{1}{}%
_{f}S^{K}(\mathbb{Q}{})$.}

\begin{proposition}
\label{p8}The group $\varprojlim\nolimits_{K\in\mathcal{K}{}}^{1}{}_{f}%
S^{K}(\mathbb{Q}{})$ is uncountable.
\end{proposition}

From (\ref{e1}) applied to the cohomology sequence of (\ref{e05}), we see
that
\[
\varprojlim\nolimits_{K\in\mathcal{K}{}}^{1}{}_{f}S^{K}(\mathbb{Q}{}%
)\overset{\cong}{\rightarrow}\varprojlim\nolimits_{K\in\mathcal{K}{}}^{1}%
{}_{f}\bar{S}^{K}(\mathbb{Q}{}),
\]
and so we compute the second group.

\begin{lemma}
\label{p9}Let $\mathcal{F}{}$ be the set of totally real subfields of
$\mathbb{C}{}$ that are finite and Galois over $\mathbb{Q}{}$. Then the
inverse system $(_{f}\bar{S}^{K})_{K\in\mathcal{K}{}}$ is equivalent with the
inverse system $((\mathbb{G}_{m})_{B(F)/\mathbb{Q}{}})_{F\in\mathcal{F}{}}$.
Thus,
\[
\varprojlim\nolimits_{K\in\mathcal{K}{}}^{1}{}_{f}\bar{S}^{K}(\mathbb{Q}%
{})\approx\varprojlim\nolimits_{F\in\mathcal{F}{}}^{1}B(F)^{\times}%
\]
(on the right, the transition maps are the norm maps).
\end{lemma}

\begin{proof}
Fix a quadratic imaginary field $k$. The fields $K\in\mathcal{K}{}$ containing
$k$ form a cofinal set. Once $k$ has been fixed, the isomorphism in Lemma
\ref{p4} becomes canonical. In particular, it is natural for the norm maps.
\end{proof}

Note that the isomorphism in the lemma depends only on the choice of $k$.

\begin{example}
\label{p10}If $\Gamma_{F/\mathbb{Q}{}}=\Gal(F/\mathbb{Q}{})$ is commutative,
$B(F)$ is a product of copies of $\mathbb{Q}{}$ indexed by the elements of
$\Gamma_{F/\mathbb{Q}{}}$. Thus, $(\mathbb{G}_{m})_{B(F)/\mathbb{Q}{}}%
\approx\mathbb{G}_{m}^{[F\colon\mathbb{Q}{}]}$.
\end{example}

\begin{example}
\label{p11}Let $F=F_{1}\cdot F_{2}$ with $F_{1},F_{2}\in\mathcal{F}{}$ and
$F_{1}\cap F_{2}=\mathbb{Q}{}$. Then%
\[
\Gal(F/\mathbb{Q}{})\cong\Gal(F_{1}/\mathbb{Q}{})\times\Gal(F_{2}/\mathbb{Q}%
{})
\]
and a conjugacy class $C$ in $\Gal(F_{1}/\mathbb{Q}{})$ is the image of the
conjugacy class $C\times\{1\}$ in $\Gal(F/\mathbb{Q})$. Therefore, $B(F_{1})$
is a direct factor of $B(F)$, and so the norm map%
\[
B(F)^{\times}\rightarrow B(F_{1})^{\times}%
\]
is surjective.
\end{example}

\begin{example}
\label{p12}Fix an odd prime $l$, and write $C_{l}$ for any cyclic group of
order $l$. Define $G$ to be the semi-direct product $N\rtimes_{\theta}Q$ of
$N=C_{l}\times C_{l}$ (generators $a,b$) with $Q=C_{l}$ (generator $c$)
relative to the homomorphism $\theta\colon C_{l}\rightarrow\Aut(C_{l}\times
C_{l})$ for which
\[
\theta(c^{i})=\left(
\begin{array}
[c]{ll}%
1 & 0\\
i & 1
\end{array}
\right)  \text{, i.e., }\theta(c^{i})(a)=ab^{i},\quad\theta(c^{i}%
)(b)=b\text{.}%
\]
Then $G$ has generators $a,b,c$, and relations%
\[
a^{l}=b^{l}=c^{l}=1,\quad ab=cac^{-1},\quad\lbrack b,a]=1=[b,c]\text{.}%
\]
All elements $\neq1$ in $G$ have order $l$, and the centre of $G$ is $\langle
b\rangle$.

The inverse image under $G\rightarrow Q$ of the conjugacy class $\{c\}$ breaks
up into $l$ conjugacy classes, namely,%
\[
\{a^{j}c,ba^{j}c,\ldots,b^{l-1}a^{j}c\},\quad0\leq j\leq l-1\text{,}%
\]
(because $a^{-1}ca=bc$). The centralizer of $a^{j}c$ in $G$ is $\langle
b,a^{j}c\rangle$, which has order $l^{2}$.

Let $E$ be an extension of $\mathbb{Q}{}$ with Galois group $G$, and let
$F_{0}=E^{N}$. Then (\ref{p10})%
\[
B(F_{0})\cong F_{0}^{\{1\}}\times F_{0}^{\{c\}}\times\cdots\times
F_{0}^{\{c^{l-1}\}},\quad F_{0}^{\{c^{i}\}}=\mathbb{Q}\text{,}%
\]
and the inverse image of $F_{0}^{\{c\}}$ in $B(E)$ under the norm map is%
\[
E^{\langle b,c\rangle}\times E^{\langle b,ac\rangle}\times\cdots\times
E^{\langle b,a^{l-1}c\rangle}.
\]

\end{example}

\begin{lemma}
\label{p13}The group $\varprojlim\nolimits_{F\in\mathcal{F}{}}^{1}%
B(F)^{\times}$ is uncountable.
\end{lemma}

\begin{proof}
According to Proposition \ref{i1}b, it suffices to show that the inverse
system $(B(F)^{\times})_{F\in\mathcal{F}{}}$ fails (ML). Examples \ref{p10}
and \ref{p11} show that we shall need to consider nonabelian Galois groups and
nonsplit extensions.

Fix an odd prime number $l$, and choose a prime number $p_{0}$ that splits
completely in $\mathbb{Q}{}[\sqrt[l]{1}]$. Then $l|p_{0}-1$, and so there is a
surjective homomorphism%
\[
(\mathbb{Z}{}/p_{0}\mathbb{Z}{})^{\times}\rightarrow C_{l}.
\]
Let $F_{0}$ be the subfield of $\mathbb{Q}{}[\sqrt[p_{0}]{1}]$ fixed by the
kernel of one such homomorphism.

To prove that $(B(F)^{\times})_{F\in\mathcal{F}{}}$ fails (ML) it suffices to
show that, for each $L\in\mathcal{F}{}$ containing $F_{0}$, there exists an
$E\in\mathcal{F}{}$ containing $F_{0}$ and such that%
\[
\Nm_{B(E)/B(F_{0})}(B(E)^{\times})\not \supset \Nm_{B(L)/B(F_{0}%
)}(B(L)^{\times})\text{;}%
\]
for then $E\cdot L\in\mathcal{F}{}$, but%
\[
\Nm_{B(E\cdot L)/B(F_{0})}(B(E\cdot L)^{\times})\not =\Nm_{B(L)/B(F_{0}%
)}(B(L)^{\times})\text{.}%
\]

Choose a prime $p_{1}$ that splits completely in $L[\sqrt[l]{1},\sqrt[l]%
{p_{0}}]$, and construct a cyclic extension $F_{1}$ of $\mathbb{Q}{}$ of
degree $l$ by choosing a surjective homomorphism $(\mathbb{Z}{}/p_{1}%
\mathbb{Z}{})^{\times}\rightarrow C_{l}$, as before. Then $F=_{\text{df}}%
F_{0}\cdot F_{1}$ has Galois group%
\[
\Gal(F/\mathbb{Q}{})\cong C_{l}\times C_{l}.
\]
Let $G$ be as in Example \ref{p12}, and consider the extension%
\[
1\rightarrow\langle b\rangle\rightarrow G\rightarrow\langle a,c\rangle
\rightarrow1.
\]
Let
\[
\alpha\colon G/\langle b\rangle\rightarrow\Gal(F/\mathbb{Q}{})
\]
be the isomorphism sending $c$ to a generator of $\Gal(F_{0}/\mathbb{Q}{})$
and $a$ to a generator of $\Gal(F_{1}/\mathbb{Q}{})$. The only primes
ramifying in $F$ are $p_{0}$ and $p_{1}$, and for $p=p_{0}$ or $p_{1}$

\begin{itemize}
\item $l$ divides $p-1$ (because both primes split in $\mathbb{Q}{}%
[\sqrt[l]{1}]$);

\item for all primes $v$ of $F$ dividing $p$, $F_{v}$ is totally ramified over
$\mathbb{Q}{}_{p}$ ($p_{i}$ is totally ramified in $F_{i}$; $p_{1}$ splits
completely in $F_{0}$ by construction; $p_{0}$ splits completely in $F_{1}$
because it becomes an $l^{\text{th}}$ power in $\mathbb{F}_{p_{1}}$).
\end{itemize}

\noindent Now an argument of Scholz and Reichardt (see Serre 1992, Theorem
2.1.3) shows that there exists a Galois extension $E$ of $\mathbb{Q}{}$
containing $F$ for which there is a commutative diagram%
\[
\begin{CD}
G @>>>G/\langle b\rangle\\
@VV{\approx}V@VV{\approx}V\\
\Gal(E/\mathbb{Q})@>>> \Gal(F/\mathbb{Q}).
\end{CD}
\]
Note that $E$, being Galois of odd degree over $\mathbb{Q}{}$, must be totally
real.%
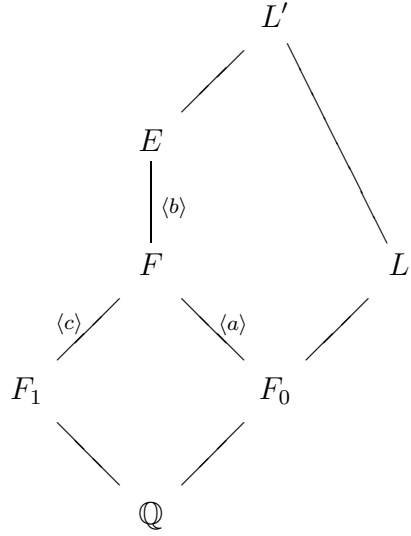
\[
\begin{diagram}
&  &  &  & L^{\prime} &  & \\
&  &  & \ldLine &  & \rdLine(2,4) & \\
&  & E &  &  &  & \\
&  & \dLine_{\langle b\rangle} &  &  &  & \\
&  & F &  &  &  & L\\
& \ldLine^{\langle c\rangle} &  & \rdLine^{\langle a\rangle} &  & \ldLine & \\
F_{1} &  &  &  & F_{0} &  & \\
& \rdLine &  & \ldLine &  &  & \\
&  & \mathbb{Q} &  &  &  &
\end{diagram}
\]

I claim that the image of $B(E)^{\times}$ in $F_{0}^{\{c\}\times}$ does not
contain the image of $B(L)^{\times}$. In order to show this, it suffices to
show that the image of $(B(E)\otimes\mathbb{Q}_{p_{1}})^{\times}$ in
$(F_{0}^{\{c\}}\otimes\mathbb{Q}{}_{p_{1}})^{\times}=\mathbb{Q}_{p_{1}%
}^{\times}$ does not contain the image of $(B(L)\otimes\mathbb{Q}{}_{p_{1}%
})^{\times}$. But, because $p_{1}$ splits completely in $L$, the second group
is $\mathbb{Q}{}_{p_{1}}^{\times}$. On the other hand, $p_{1}$ is totally
ramified in each field $E^{\langle b,a^{j}c\rangle}$ (because $E^{\langle
b,a^{j}c\rangle}\cdot F_{0}=F$, $(p_{1})=\mathcal{\mathfrak{p}}_{1}^{l}%
\cdots\mathfrak{p}{}_{l}^{l}$ in $F$, and $p_{1}$ splits completely in $F_{0}%
$), and so, for each $v|p$, $E_{v}^{\langle b,a^{j}c\rangle}$ is the (unique)
tamely ramified cyclic extension of $\mathbb{Q}{}_{p_{1}}$ of degree $l$.
Thus, the image of $\prod_{v}E_{v}^{\langle b,a^{j}c\rangle\times}$ in
$\mathbb{Q}{}_{p_{1}}^{\times}$ does not contain $\mathbb{Z}{}_{p_{1}}%
^{\times}$.
\end{proof}

\subsubsection{Characterizing the period torsor}

\begin{theorem}
\label{p14}Let $K$ be a CM-subfield of $\mathbb{C}$, Galois over $\mathbb{Q}%
{}$.

\begin{enumerate}
\item If $[K\colon\mathbb{Q}{}]<\infty$, then $\mathsf{P}(G_{\text{B}}%
^{\CM,K}\rightarrow\varGamma;P^{\Art})$ contains exactly one isomorphism
class, which is represented by $P^{\CM,K}\rightarrow P^{\Art}$.

\item If $K=\mathbb{Q}{}^{\text{cm}}$, then $\mathsf{P}(G_{\text{B}}%
^{\CM,K}\rightarrow\varGamma;P^{\Art})$ contains uncountably many isomorphism classes.
\end{enumerate}
\end{theorem}

\begin{proof}
According to Proposition \ref{p00}, the isomorphism classes in $\mathsf{P}%
(G_{\text{B}}^{\CM,K}\rightarrow\varGamma;P^{\Art})$ are classified by
$H^{1}(\mathbb{Q},{}_{f}S^{K})$. Therefore, (a) and (b) follow respectively
from Proposition \ref{p5} and Proposition \ref{p01}.
\end{proof}

\begin{remark}
\label{p15}In fact, we have shown that the uncountable group $\varprojlim
_{K}^{1}{}_{f}S^{K}(\mathbb{Q}{})$ acts simply transitively on the set of
isomorphism classes of objects in $\mathsf{P}(G_{\text{B}}^{\CM}%
\rightarrow\varGamma;P^{\Art})$. To make this explicit, let $K_{0}$ be a
quadratic imaginary field, and let%
\[
K_{0}\subset\cdots\subset K_{n}\subset\cdots\subset\mathbb{Q}{}^{\text{cm}}%
\]
be a sequence of CM-fields with union $\mathbb{Q}{}^{\text{cm}}$. Let
$S_{n}=S^{K_{n}}$ and $P_{n}^{\CM}=P^{\CM,K_{n}}$. From an element
$s=(s_{n})_{n\in\mathbb{N}}$ of $\prod{}_{f}S_{n}(\mathbb{Q}{})$ we obtain a
$v$-morphism of torsors $P^{\CM}(s)\rightarrow P^{\Art}$ by modifying the
transition maps in $P^{\CM}$: define%
\[
P^{\CM}(s)=\varprojlim_{n\in\mathbb{N}{}}(P_{n}^{\CM},v_{n}\circ
s_{n})\text{.}%
\]
Then, the isomorphism class of $P^{\CM}(s)\rightarrow P^{\Art}$ depends only
on the class of $(s_{n})$ in $\varprojlim_{n}^{1}{}_{f}S_{n}(\mathbb{Q}{})$,
distinct classes in $\varprojlim_{n}^{1}{}_{f}S_{n}(\mathbb{Q}{})$ give
nonisomorphic objects in $\mathsf{P}(G_{B}^{\CM}\rightarrow\varGamma;P^{\Art}%
)$, and every object in $\mathsf{P}(G_{B}^{\CM}\rightarrow\varGamma;P^{\Art})$
is isomorphic to $P^{\CM}(s)\rightarrow P^{\Art}$ for some $s\in\prod{}%
_{f}S_{n}(\mathbb{Q}{})$.
\end{remark}

\begin{remark}
\label{p16}There remains the problem of characterizing the isomorphism class
of $P^{\CM}\rightarrow P^{\Art}$. One may hope that it is uniquely determined
by its isomorphism classes in $\mathsf{P}((G_{\text{B}}^{\CM})_{/\mathbb{Q}%
{}_{l}}\rightarrow\varGamma_{/\mathbb{Q}{}_{l}},P_{/\mathbb{Q}{}_{l}}^{\Art})$
for $l=2,3,\ldots,\infty$, i.e., that the kernel of%
\[
H^{1}(\mathbb{Q}{},{}_{f}S)\rightarrow%
%TCIMACRO{\tprod _{l}}%
%BeginExpansion
{\textstyle\prod_{l}}
%EndExpansion
H^{1}(\mathbb{Q}{}_{l},{}_{f}S)
\]
is zero, but the calculations I have made in this direction do not look
promising. Note that ${}_{f}S_{/\mathbb{R}{}}$ is a countable product of
copies of $\mathbb{G}_{m}$, and so $H^{1}(\mathbb{R}{},{}_{f}S)=0$.
\end{remark}

\newpage

\section{Periods of abelian varieties}

For a class of abelian varieties $\mathcal{A}{}{}$ over $\mathbb{Q}{}$,
$\langle\mathcal{A}{}\rangle^{\otimes}$ denotes the category of motives based
on the abelian varieties in $\mathcal{A}{}$, the projective spaces, and the
zero-dimensional varieties. Let
\begin{align*}
G_{\text{B}}^{\mathcal{A}{}{}}  &  =\mathcal{A}{}ut^{\otimes}(\omega
_{\text{B}}|\langle\mathcal{A}{}\rangle^{\otimes}),\\
G_{\text{dR}}^{\mathcal{A}{}}  &  =\mathcal{A}{}ut^{\otimes}(\omega
_{\text{dR}}|\langle\mathcal{A}{}\rangle^{\otimes}),\\
P^{\mathcal{A}{}}  &  =\mathcal{I}{}som^{\otimes}(\omega_{\text{B}}%
,\omega_{\text{dR}}).
\end{align*}
When $\mathcal{A}{}=\{A\}$, we write $A$ for $\mathcal{A}{}$.

\begin{plain}
\label{period1}The inclusion of the Artin motives into $\langle\mathcal{A}%
{}\rangle^{\otimes}$ defines a homomorphism%
\[
G_{\text{B}}^{\mathcal{A}{}}\rightarrow\varGamma\text{.}%
\]
This homomorphism is surjective, and its kernel is the identity component
$(G_{\text{B}}^{\mathcal{A}{}})^{\circ}$ of $G_{\text{B}}^{\mathcal{A}{}}$
(Deligne and Milne 1982, 6.23)\footnote{\smallskip Which applies because of
Deligne's theorem (Deligne 1982a) that all Hodge classes on abelian varieties
are absolutely Hodge.}. In particular, for a single abelian variety $A$, there
is an exact sequence%
\[
1\rightarrow\MT(A)\rightarrow G_{\text{B}}^{A}\rightarrow\varGamma\rightarrow
1\text{.}%
\]

\end{plain}

\begin{plain}
\label{period2}For a single abelian variety $A$, there is a unique
homomorphism%
\[
S\rightarrow\MT(A)_{/\mathbb{R}{}}^{\text{ab}}%
\]
sending $h_{\text{can}}$ onto $(h^{A})^{\text{ab}}$ where, as usual,
$h^{A}\colon\mathbb{S}{}\rightarrow\MT(A)_{/\mathbb{R}{}}$ is the homomorphism
defining the Hodge structure on $H^{\ast}(A,\mathbb{Q}{})$. This homomorphism
is surjective, and it factors through $S^{K}$ if and only if $K$ contains the
reflex field of $(\MT(A),h^{A})^{\text{ab}}$.
\end{plain}

\begin{plain}
\label{period3}On combining the last two statements, we see that if $A$ is an
abelian variety over $\mathbb{Q}{}$ such that $S^{K}\rightarrow
\MT(A)^{\text{ab}}$ is an isomorphism for some $K$, then%
\[
1\rightarrow\MT(A)^{\text{der}}\rightarrow G_{\text{B}}^{A}\rightarrow
G_{\text{B}}^{\CM,K}\rightarrow1
\]
is exact.
\end{plain}

\subsection{Some Hasse principles}

We first prove an elementary structure theorem.

\begin{lemma}
\label{hp}Every semisimple group $H$ over $\mathbb{Q}{}$ such that
$H_{/\mathbb{Q}{}^{\text{al}}}$ is a product of simple groups is isomorphic to
a product of groups $H_{i}$ of the form $H_{i}=\Res_{F_{i}/\mathbb{Q}{}}N_{i}$
with $F_{i}$ a number field and $N_{i}$ an absolutely simple group over
$F_{i}$.
\end{lemma}

\begin{proof}
To give a semisimple group over $\mathbb{Q}{}$ is the same as to give a
semisimple group $H$ over $\mathbb{Q}{}^{\text{al}}$ together with a descent
datum $(\alpha_{\sigma})_{\sigma\in\Gamma}$. Here $\alpha_{\sigma}$ is an
isomorphism $\sigma H\rightarrow H$, $\alpha_{\sigma}\circ\sigma\alpha_{\tau
}=\alpha_{\sigma\tau}$ for all $\sigma$ and $\tau$, and there is a continuity condition.

First consider a pair $\left(  H,(\alpha_{\sigma})_{\sigma\in\Gamma}\right)  $
with $H$ an adjoint group. Write $H=%
%TCIMACRO{\tprod _{i\in I}}%
%BeginExpansion
{\textstyle\prod_{i\in I}}
%EndExpansion
H_{i}$ as a product of simple groups. For each $\sigma\in\Gamma$, there is a
permutation (also denoted $\sigma$) of $I$ such that $\alpha_{\sigma}$ is a
product of isomorphisms
\[
\alpha_{\sigma}(i)\colon\sigma H_{i}\rightarrow H_{\sigma i}\text{.}%
\]
Let $J$ be an orbit of $\Gamma$ in $I$, let $j\in J$, and let $\Gamma
_{j}=\{\sigma\in\Gamma\mid\sigma j=j\}$. Then $(\alpha_{\sigma}(j))_{\sigma
\in\Gamma_{j}}$ is a descent datum on $H_{j}$, and $\left(
%TCIMACRO{\tprod _{i\in J}}%
%BeginExpansion
{\textstyle\prod_{i\in J}}
%EndExpansion
\alpha_{\sigma}(i)\right)  _{\sigma\in\Gamma}$ is a descent datum on $%
%TCIMACRO{\tprod _{i\in J}}%
%BeginExpansion
{\textstyle\prod_{i\in J}}
%EndExpansion
H_{i}$. The first defines a model $N_{j}$ of $H_{j}$ over $F_{j}=_{\text{df}%
}\mathbb{Q}{}^{\Gamma_{j}}$, which is absolutely simple, and the second
defines a model $M_{J}$ of $%
%TCIMACRO{\tprod _{i\in J}}%
%BeginExpansion
{\textstyle\prod_{i\in J}}
%EndExpansion
H_{i}$ over $\mathbb{Q}{}$, which is isomorphic to $\Res_{F_{j}/\mathbb{Q}{}%
}N_{j}$. Now $%
%TCIMACRO{\tprod _{J\in\Gamma\backslash I}}%
%BeginExpansion
{\textstyle\prod_{J\in\Gamma\backslash I}}
%EndExpansion
M_{J}$ is a semisimple group over $\mathbb{Q}{}$ giving rise to $(H,(\alpha
_{\sigma})_{\sigma\in\Gamma})$ over $\mathbb{Q}{}^{\text{al}}$.

Next consider a pair $\left(  H,(\alpha_{\sigma})_{\sigma\in\Gamma}\right)  $
with $H$ a product $H=%
%TCIMACRO{\tprod _{i\in I}}%
%BeginExpansion
{\textstyle\prod_{i\in I}}
%EndExpansion
H_{i}$ of simple groups. Then $(\alpha_{\sigma})_{\sigma}$ defines a descent
datum $(\alpha_{\sigma}^{\text{ad}})_{\sigma}$ on $H^{\text{ad}}$, and, as
above, $\alpha_{\sigma}^{\text{ad}}$ is a product of isomorphisms
$\alpha_{\sigma}^{\text{ad}}(i)\colon\sigma H_{i}^{\text{ad}}\rightarrow
H_{\sigma i}^{\text{ad}}$. Consider%
\[
\begin{CD}
\sigma H_{i} @>{\alpha_{\sigma}(i)}>> H_{\sigma i}\\
@VVV @VVV\\
\sigma H_{i}^{\text{ad}} @>{\alpha_{\sigma}^{\text{ad}}(i)}>> H_{\sigma
i}^{\text{ad}}.
\end{CD}
\]
Here $\alpha_{\sigma}(i)$ is the composite%
\[
\sigma H_{i}\hookrightarrow\sigma H\overset{\alpha_{\sigma}}{\rightarrow
}H\xr{\textrm{project}}H_{i}\text{.}%
\]
Because the diagram commutes, $\alpha_{\sigma}$ and $%
%TCIMACRO{\tprod _{i\in I}}%
%BeginExpansion
{\textstyle\prod_{i\in I}}
%EndExpansion
\alpha_{\sigma}(i)$ differ by a map from $\sigma H$ into the centre of $H$,
which must be trivial because $\sigma H$ is connected. Thus, $\alpha_{\sigma}=%
%TCIMACRO{\tprod _{i\in I}}%
%BeginExpansion
{\textstyle\prod_{i\in I}}
%EndExpansion
\alpha_{\sigma}(i)$, and the same argument as in the preceding paragraph
completes the proof of the lemma.
\end{proof}

We define the \emph{index (of connectivity) }of a semisimple algebraic group
$H$ to be the degree of the universal covering $\tilde{H}\rightarrow H$. Thus,
for an isogeny $a\colon H^{\prime}\rightarrow H$ of semisimple groups,%
\[
\text{index}(H^{\prime})\leq\text{index}(H)
\]
with equality if and only if $a$ is an isomorphism.

Consider the following condition on a semisimple algebraic group $H$ over
$\mathbb{Q}{}$:

\begin{quotation}
(*) $H_{/\mathbb{Q}{}^{\text{al}}}$ is a product of simple groups of index $1$
or $2$.
\end{quotation}

\begin{proposition}
\label{g03}For any semisimple group $H$ over $\mathbb{Q}{}$ satisfying (*),
the map%
\[
H^{1}(\mathbb{Q}{},H)\rightarrow%
%TCIMACRO{\tprod _{l=2,\ldots,\infty}}%
%BeginExpansion
{\textstyle\prod_{l=2,\ldots,\infty}}
%EndExpansion
H^{1}(\mathbb{Q}_{l},H)
\]
is injective.
\end{proposition}

\begin{proof}
According to (\ref{hp}), $H\approx%
%TCIMACRO{\tprod }%
%BeginExpansion
{\textstyle\prod}
%EndExpansion
\Res_{F_{i}/\mathbb{Q}{}}N_{i}$ with each $N_{i}$ absolutely simple of index
$1$ or $2$, and%
\[
H^{1}(\mathbb{Q}{},H)\approx%
%TCIMACRO{\tprod _{i}}%
%BeginExpansion
{\textstyle\prod_{i}}
%EndExpansion
H^{1}(\mathbb{Q}{},H_{i})\approx%
%TCIMACRO{\tprod _{i}}%
%BeginExpansion
{\textstyle\prod_{i}}
%EndExpansion
H^{1}(F_{i},N_{i}).
\]
Therefore, the proposition follows from the next two lemmas.
\end{proof}

\begin{lemma}
\label{g02}For any simply connected semisimple group $H$ over a number field
$F{}$, the map%
\[
H^{1}(F,H)\rightarrow%
%TCIMACRO{\tprod _{v\text{ real}}}%
%BeginExpansion
{\textstyle\prod_{v\text{ real}}}
%EndExpansion
H^{1}(F_{v},H)
\]
is bijective.
\end{lemma}

\begin{proof}
This is the theorem of Kneser, Harder, and Chernousov --- see Platonov and
Rapinchuk 1994, Theorem 6.6, p286.
\end{proof}

\begin{lemma}
\label{g1}For any semisimple group $H$ of index $2$ over a number field $F$,
\[
H^{1}(F,H)\rightarrow%
%TCIMACRO{\tprod _{v}}%
%BeginExpansion
{\textstyle\prod_{v}}
%EndExpansion
H^{1}(F_{v},H)\quad(v\text{ runs over all primes of }F)
\]
is injective.
\end{lemma}

\begin{proof}
Platonov and Rapinchuk (1994, Remark p337) note that the map has trivial
kernel. The lemma can now be proved by a twisting argument, because any form
of a semisimple group of index $2$ again has index $2$.
\end{proof}

\begin{proposition}
\label{g04}Let $H$ be a semisimple group over $\mathbb{Q}{}$ satisfying (*).
If two cohomology classes in $H^{1}(\mathbb{Q}{},\tilde{H})$ have the same
image in $H^{1}(\mathbb{R}{},H)$, then they have the same image in
$H^{1}(\mathbb{Q}{},H)$. In other words, $H^{1}(\mathbb{Q}{}{},H)\rightarrow
H^{1}(\mathbb{R}{},H)$ is injective on the image of $H^{1}(\mathbb{Q}{}%
,\tilde{H})\rightarrow H^{1}(\mathbb{Q}{},H)$.
\end{proposition}

\begin{proof}
In the diagram%
\[
\begin{CD}
H^{1}(\mathbb{Q},Z) @>>> H^{1}(\mathbb{Q},\tilde{H}) @>>> H^{1}(\mathbb{Q},H)\\
@VV{\text{surjective}}V @VV{\text{injective}}V @VVV \\
H^{1}(\mathbb{R},Z) @>>>H^{1}(\mathbb{R},\tilde{H}) @>>>
H^{1}(\mathbb{R},H),
\end{CD}
\]
the rows are exact in the sense that the fibres of the second map are the
orbits of the natural action of the first group on the middle set. The first
vertical arrow is surjective, and the second is injective. A diagram chase
completes the proof.
\end{proof}

\begin{remark}
\label{g05}The proof of Proposition \ref{g04} shows that for any central
extension%
\[
1\rightarrow Z\rightarrow G^{\prime}\rightarrow G\rightarrow1
\]
of algebraic groups, the map $H^{1}(\mathbb{Q}{},G)\rightarrow H^{1}%
(\mathbb{R}{},G)$ is injective on the image of $H^{1}(\mathbb{Q}{},G^{\prime
})$ in $H^{1}(\mathbb{Q}{},G)$ provided

\begin{enumerate}
\item the map $H^{1}(\mathbb{Q}{},Z)\rightarrow H^{1}(\mathbb{R}{},Z)$ is surjective;

\item the map $H^{1}(\mathbb{Q}{},G^{\prime})\rightarrow H^{1}(\mathbb{R}%
{},G^{\prime})$ is injective.
\end{enumerate}
\end{remark}

\begin{proposition}
\label{g23}Let $G$ be connected reductive group over $\mathbb{Q}{}$ such that
$G^{\text{der}}$ is simply connected.

\begin{enumerate}
\item If $H^{1}(\mathbb{Q}{},G^{\text{ab}})\rightarrow%
%TCIMACRO{\tprod _{l}}%
%BeginExpansion
{\textstyle\prod_{l}}
%EndExpansion
H^{1}(\mathbb{Q}_{l}{},G^{\text{ab}})$ is injective, then so also is
$H^{1}(\mathbb{Q}{},G)\rightarrow%
%TCIMACRO{\tprod _{l}}%
%BeginExpansion
{\textstyle\prod_{l}}
%EndExpansion
H^{1}(\mathbb{Q}_{l}{},G)$ (product over $l=2,3,5,\ldots,\infty$).

\item If $H^{1}(\mathbb{Q}{},G^{\text{ab}})\rightarrow H^{1}(\mathbb{R}%
{},G^{\text{ab}})$ is injective, then so also is $H^{1}(\mathbb{Q}%
{},G)\rightarrow H^{1}(\mathbb{R}{},G).$
\end{enumerate}
\end{proposition}

\begin{proof}
(a) See Deligne 1971, 5.12. (Recall that, because $G^{\text{der}}$ is simply
connected, $H^{1}(\mathbb{Q}{}_{l},G^{\text{der}})=0$ for $l\neq\infty$. Using
this, we obtain a commutative diagram with exact rows%
\[
\begin{CD}
@.G^{\text{ab}}(\mathbb{Q}{})@>>>H^{1}(\mathbb{Q},G^{\text{der}})@>>> H^{1}(\mathbb{Q}{},G) @>>>H^{1}(\mathbb{Q}{},G^{\text{ab}})\\
@.@VVV@VVV@VVV@VVV\\
G(\mathbb{R})@>>> G^{\text{ab}}(\mathbb{R}{})@>>>H^{1}(\mathbb{R}{},G^{\text{der}}) & \rightarrow &
%TCIMACRO{\tprod _{l}}%
%BeginExpansion
{\textstyle\prod_{l}}
%EndExpansion
H^{1}(\mathbb{Q}_{l}{},G) @>>>
%TCIMACRO{\tprod _{l}}%
%BeginExpansion
{\textstyle\prod_{l}}
%EndExpansion
H^{1}(\mathbb{Q}_{l},G^{\text{ab}})\text{.}%
\end{CD}
\]
The image of $G(\mathbb{R}{})$ in $G^{\text{ab}}(\mathbb{R}{})$ contains its
identity component, and the real approximation theorem shows that
$G^{\text{ab}}(\mathbb{R}{})$ maps onto $\pi_{0}(G^{\text{ab}}(\mathbb{R}{}%
))$. Now a diagram chase shows that any element of $H^{1}(\mathbb{Q}{},G)$
that is locally zero is zero. To show that two elements $c$ and $c^{\prime}$
of $H^{1}(\mathbb{Q}{},G)$ are equal when they are locally, choose a torsor
representing $c$, and twist the groups by it.)

(b) Similar, but easier.
\end{proof}

\subsection{Characterizing $P^{\AV}\rightarrow P^{\CM}$}

Let $\AV(\mathbb{Q}{})=\langle\mathcal{A}{}\rangle^{\otimes}$ with
$\mathcal{A}{}$ the class of all abelian varieties over $\mathbb{Q}{}$.

\begin{theorem}
\label{g01}Let $A$ be an abelian variety over $\mathbb{Q}{}$. If
$\MT(A)^{\text{der}}$ satisfies (*) and\footnote{We don't really need to
assume $\MT(A)^{\text{ab}}\cong S^{K}$ --- it only makes the statement a
little more pleasant.} $\MT(A)^{\text{ab}}\cong S^{K}$ for some $K$, then the
isomorphism class of $P^{A}\rightarrow P^{\CM,K}$ is uniquely determined by
its classes over $\mathbb{Q}{}_{l}$ ($l=2,3,\ldots,\infty$).
\end{theorem}

\begin{proof}
According to Proposition \ref{t4}, the isomorphism classes in $\mathsf{P}%
(G_{\text{B}}^{A}\rightarrow G_{\text{B}}^{\CM,K};P^{\CM,K})$ are classified
by $H^{1}(\mathbb{Q}{},{}H)$ where $H$ is the twist of $\Ker(G_{\text{B}}%
^{A}\rightarrow G_{\text{B}}^{\CM,K})$ by the period torsor. But according to
(\ref{period3}), this kernel is $\MT(A)^{\text{der}}$, and so the theorem
follows from Proposition \ref{g03}.
\end{proof}

The next result shows that the abelian varieties over $\mathbb{Q}{}$
satisfying the conditions of Theorem \ref{g01} are cofinal among all abelian
varieties over $\mathbb{Q}{}$ for the relation $A\prec B$ if $hA$ is
isomorphic to an object of $\langle B\rangle^{\otimes}$.

\begin{proposition}
\label{g6}Let $A$ be an abelian variety over $\mathbb{Q}{}$. For any
sufficiently large CM-subfield $K$ of $\mathbb{C}{}$, there exists an abelian
variety $B$ over $\mathbb{Q}{}$ such that $\MT(A\times B)^{\text{der}}$
satisfies (*) and $\MT(A\times B)^{\text{ab}}\cong S^{K}$.
\end{proposition}

This will be proved in the next subsubsection.

\begin{theorem}
\label{g7}The isomorphism class of $P^{\AV}\rightarrow P^{\CM}$ is uniquely
determined by its classes over $\mathbb{Q}{}_{l}$ ($l=2,3,\ldots,\infty$).
\end{theorem}

\begin{proof}
According to Proposition \ref{t4}, the isomorphism classes in $\mathsf{P}%
(G_{\text{B}}^{\AV}\rightarrow G_{\text{B}}^{\CM};P^{\CM})$ are classified by
$H^{1}(\mathbb{Q}{},H)$ where $H$ is the twist of $H^{\prime}=\Ker(G_{\text{B}%
}^{\AV}\rightarrow G_{\text{B}}^{\CM})$ by the period torsor. Proposition
\ref{g6} implies that
\[
H^{\prime}\cong\varprojlim\MT(A)^{\text{der}}%
\]
where $A$ runs over the abelian varieties satisfying the hypothesis of Theorem
\ref{g01}. Therefore $H$ is a countable product of simple algebraic groups
satisfying (*), and so we can apply Corollary \ref{c3} and Proposition
\ref{g03}.
\end{proof}

\subsubsection{Proof of Proposition \ref{g6}}

We shall say that an algebraic group $H$ over a field $k$ of characteristic
zero is a \emph{special orthogonal group} if there exists vector space $V$
over $k$ of dimension $n>8$ and a nondegenerate quadratic form $q$ on $V$ such
that $H\approx\SO(V,q)$. Such an $H$ is a connected and simple, and there are
isogenies of degree $2$%
\[
\tilde{H}\rightarrow H\rightarrow H^{\text{ad}}%
\]
with $\tilde{H}$ the simply connected covering group of $H$ (a spinor group).
Let $Z=Z\tilde{H}$. When $n$ is even $Z(k^{\text{al}})\approx C_{2}\times
C_{2}$.

\begin{plain}
\label{g06}Recall from Deligne 1979a the following statements.

\begin{enumerate}
\item Let $A$ be an abelian variety over $\mathbb{C}{}$. Each $\mathbb{Q}{}%
$-simple factor $H$ of $\MT(A)^{\text{ad}}$ is of type $A$, $B$, $C$,
$D^{\mathbb{R}{}}$, or $D^{\mathbb{H}{}}$ (notations as in ibid. 2.3.8).

\item Let $A$ be an abelian variety over $\mathbb{C}{}$, and let $G=\MT(A)$.
The homomorphism $h\colon\mathbb{S}{}\rightarrow G_{/\mathbb{R}{}}$ defined by
the Hodge structure on $H^{\ast}(A,\mathbb{Q}{})$ satisfies the following conditions:

\begin{enumerate}
\item the Hodge structure on the Lie algebra $\mathfrak{g}{}$ of $G$ defined
by $\Ad\circ h\colon\mathbb{S}{}\rightarrow\GL(\mathfrak{g}{})$ is of type
$\{(1,-1),(0,0),(-1,1)\}$;

\item $\ad h(i)$ is a Cartan involution on $G^{\text{ad}}$;

\item $h$ generates $G$ (i.e., there is no proper closed subgroup $G^{\prime
}\subset G$ such that $h(\mathbb{S}{})\subset G_{/\mathbb{R}{}}^{\prime}$).
\end{enumerate}

\item Let $H$ be a simple adjoint group over $\mathbb{Q}{}$, and let
$h\colon\mathbb{S}/\mathbb{G}_{m}\rightarrow H_{/\mathbb{R}{}}$ be a
homomorphism satisfying the conditions (b) i), ii), iii).

\begin{enumerate}
\item If $H$ is of type $A$, $B$, $C$, or $D^{\mathbb{R}{}}$, there exists an
abelian variety $A$ over $\mathbb{C}{}$ such that $\MT(A)^{\text{der}}$ is
simply connected and $(\MT(A),h^{A})^{\text{ad}}\approx(H,h)$ (apply ibid. 2.3.10).

\item \label{d1}Suppose $H$ is of type $D^{\mathbb{H}{}}$. Then
$H=\Res_{F/\mathbb{Q}{}}N$ for some absolutely simple group $N$ over a totally
real field $F$, and we let $H^{\prime}=\Res_{F/\mathbb{Q}{}}N^{\prime}$ where
$N^{\prime}\rightarrow N$ is the double covering of $N$ that is an inner form
of a special orthogonal group (ibid. 2.3.8). Any homomorphism $\tilde
{H}\rightarrow\MT(A)$ sending $h$ to $(h^{A})^{\text{ad}}$ factors through
$H^{\prime}$ (ibid. 1.3.10). There exists an abelian variety $A$ over
$\mathbb{C}{}$ with $\MT(A)^{\text{der}}\approx H^{\prime}$ and $(\MT(A),h^{A}%
)^{\text{ad}}\approx(H,h)$ (apply ibid. 2.3.10).
\end{enumerate}
\end{enumerate}
\end{plain}

We shall need the following condition on a semisimple group $H$ over
$\mathbb{Q}{}$:

\begin{quotation}
(**) $H$ is a product of simple groups; a simple factor of $H_{\mathbb{R}{}}$
has index $2$ if it is of type $D^{\mathbb{H}{}}$ and index $1$ otherwise.
\end{quotation}

\begin{definition}
\label{z0}An abelian variety defined over a subfield of $\mathbb{C}{}$ is
\emph{maximal }if $\MT(A)^{\text{der}}$ satisfies (**).
\end{definition}

\begin{lemma}
\label{z1}A finite product of maximal abelian varieties is maximal.
\end{lemma}

\begin{proof}
Let $A=%
%TCIMACRO{\tprod _{j}}%
%BeginExpansion
{\textstyle\prod_{j}}
%EndExpansion
A_{j}$ where each $A_{j}$ is maximal. The canonical map $\MT(A)\rightarrow%
%TCIMACRO{\tprod }%
%BeginExpansion
{\textstyle\prod}
%EndExpansion
\MT(A_{j})$ is injective and its composite with any projection $%
%TCIMACRO{\tprod }%
%BeginExpansion
{\textstyle\prod}
%EndExpansion
\MT(A_{j})\rightarrow\MT(A_{j})$ is surjective (Deligne 1982a, \S 3). For a
judicious choice of simple factors $H_{i}$ of $\prod\MT(A_{j})^{\text{der}}$,
the homomorphism $\MT(A)^{\text{ad}}\rightarrow\prod H_{i}^{\text{ad}}$ will
be an isomorphism. Consider%
\[%
%TCIMACRO{\tprod }%
%BeginExpansion
{\textstyle\prod}
%EndExpansion
\tilde{H}_{i}\rightarrow\MT(A)^{\text{der}}\rightarrow%
%TCIMACRO{\tprod }%
%BeginExpansion
{\textstyle\prod}
%EndExpansion
H_{i}\text{.}%
\]
Then $%
%TCIMACRO{\tprod }%
%BeginExpansion
{\textstyle\prod}
%EndExpansion
\tilde{H}_{i}\rightarrow\MT(A)^{\text{der}}$ factors through $%
%TCIMACRO{\tprod }%
%BeginExpansion
{\textstyle\prod}
%EndExpansion
H_{i}^{\prime}$ where $H_{i}^{\prime}$ is as in (\ref{g06}\ref{d1}) when
$H_{i}$ is of type $D^{\mathbb{H}{}}$ and equals $\tilde{H}_{i}$ otherwise.
Consider%
\[%
%TCIMACRO{\tprod }%
%BeginExpansion
{\textstyle\prod}
%EndExpansion
H_{i}^{\prime}\twoheadrightarrow\MT(A)^{\text{der}}\twoheadrightarrow%
%TCIMACRO{\tprod }%
%BeginExpansion
{\textstyle\prod}
%EndExpansion
H_{i}\text{.}%
\]
As $%
%TCIMACRO{\tprod }%
%BeginExpansion
{\textstyle\prod}
%EndExpansion
H_{i}^{\prime}$ and $%
%TCIMACRO{\tprod }%
%BeginExpansion
{\textstyle\prod}
%EndExpansion
H_{i}$ have the same index, the composite is an isomorphism, and so
$\MT(A)^{\text{der}}\rightarrow%
%TCIMACRO{\tprod }%
%BeginExpansion
{\textstyle\prod}
%EndExpansion
H_{i}$ is an isomorphism.
\end{proof}

\begin{remark}
\label{z1c}Let $A$ and $B$ be abelian varieties such that $\MT(A\times
B)^{\text{ad}}\cong\MT(B)^{\text{ad}}$. If $B$ is maximal, then so also is
$A\times B$. This can be proved the same way as Lemma \ref{z1} --- one only
has to observe that, because of the condition on the adjoint groups, the
$H_{i}$ in the proof can be chosen to be factors of $\MT(B)^{\text{der}}$.
\end{remark}

\begin{lemma}
\label{z2}For any abelian variety $A$ over $\mathbb{C}{}$, there exists an
abelian variety $B$ over $\mathbb{C}{}$ such that $(\MT(B),h^{B})^{\text{ad}%
}\approx$ ($\MT(A),h^{A})^{\text{ad}}$ and $B$ is maximal.
\end{lemma}

\begin{proof}
According to (\ref{g06}), for each simple factor $(H,h)$ of $(\MT(A),h^{A}%
)^{\text{ad}}$, there exists an abelian variety $B(H)$ such that
$(\MT(B(H)),h^{B(H)})^{\text{ad}}\approx(H,h)$ and $\MT(B(H))^{\text{der}}$ is
the covering $H^{\prime}$ of $H$ in (\ref{g06}\ref{d1}) if $H$ is of type
$D^{\mathbb{H}{}}$ and is simply connected otherwise. Take $B=%
%TCIMACRO{\tprod }%
%BeginExpansion
{\textstyle\prod}
%EndExpansion
B(H)$. Because each $B(H)$ is maximal, so is $B$ (\ref{z1}).
\end{proof}

Recall that, when $A\prec B$, there is a canonical surjection $c_{A,B}%
\colon\MT(B)\rightarrow\MT(A)$ carrying $h^{B}$ to $h^{A}$.

\begin{lemma}
\label{z3}An abelian variety $A$ over $\mathbb{C}{}$ is maximal if and only
if, for all $B$ such that $A\prec B$ and $c_{A,B}^{\text{ad}}\colon
\MT(B)^{\text{ad}}\rightarrow\MT(A)^{\text{ad}}$ is an isomorphism,
$c_{A,B}^{\text{der}}\colon\MT(B)^{\text{der}}\rightarrow\MT(A)^{\text{der}}$
is an isomorphism.
\end{lemma}

\begin{proof}
Let $B^{\prime}$ be as in (\ref{z2}), and let $B=A\times B^{\prime}$. Then
$A\prec B$ and $\MT(B)^{\text{ad}}\rightarrow\MT(A)^{\text{ad}}$ is an
isomorphism, but $\MT(B)^{\text{der}}\rightarrow\MT(A)^{\text{der}}$ is an
isomorphism only if $A$ is maximal. This proves the \textquotedblleft
if\textquotedblright.

Let $A$ be maximal. Then it is clear from (\ref{g06}c) that it satisfies the condition.
\end{proof}

\begin{lemma}
\label{z4}Let $A$ be an abelian variety over $\mathbb{C}$, and let $\tau$ be
an automorphism of $\mathbb{C}{}$. Then $A$ is maximal if and only if $\tau A$
is maximal.
\end{lemma}

\begin{proof}
Recall that $\MT(A)=\mathcal{A}{}ut^{\otimes}(\omega_{B})$ and $\MT(\tau
A)=\mathcal{A}{}ut^{\otimes}(\omega_{\tau})$ where $\omega_{\tau}$ is the
fibre functor on $\AV(\mathbb{C}{})$ sending $hX$ to $H^{\ast}(\tau
X,\mathbb{Q}{})$. Therefore, from the theory of Tannakian categories,
$\MT(\tau A)\cong{}^{P}\MT(A)$ where $P=\mathcal{I}{}som(\omega_{\text{B}%
},\omega_{\tau})$.

If $A$ is not maximal, then there exists an abelian variety $B$ such that
$A\prec B$ and $\MT(A)^{\text{ad}}\cong\MT(B)^{\text{ad}}$ but
$\MT(A)^{\text{der}}\not \cong \MT(B)^{\text{der}}$. Clearly $\tau B$ has the
same properties relative to $\tau A$, which proves that $\tau A$ is not maximal.
\end{proof}

\begin{lemma}
\label{z5}Let $A_{1}$ and $A_{2}$ be abelian varieties over $\mathbb{C}{}$
such that $(\MT(A_{1}),h^{A_{1}})^{\text{ad}}\approx(\MT(A_{2}),h^{A_{2}%
})^{\text{ad}}$. If $A_{1}$ can be defined over $\mathbb{Q}^{\text{al}}$, then
so also can $A_{2}$.
\end{lemma}

\begin{proof}
See Blasius and Borovoi 1999, 3.3.
\end{proof}

\begin{lemma}
\label{z6}For any CM-subfield $K$ of $\mathbb{C}$ of finite degree over
$\mathbb{Q}{}$, there exists an abelian variety $A$ over $\mathbb{C}$ such
that the canonical map $S^{K}\rightarrow\MT(A)$ is an isomorphism. If $K$ is
Galois over $\mathbb{Q}{}$, then $A$ may be chosen to be defined over
$\mathbb{Q}{}$.
\end{lemma}

\begin{proof}
\footnote{For slightly weaker results, see Milne 1990, I 4.6; Wei 1994, 1.5.1;
Borovoi and Blasius 1999, 3.5.}First recall that the CM-types on $K$ generate
$X^{\ast}(S^{K})$: in fact, if $\psi=\tau_{1}+\cdots+\tau_{g}$ is one CM-type
on $K$, then the CM-types $\psi_{i}=\tau_{i}+%
%TCIMACRO{\tsum _{j\neq i}}%
%BeginExpansion
{\textstyle\sum_{j\neq i}}
%EndExpansion
\iota\tau_{j}$ ($i=1,\ldots,g$) and $\bar{\phi}=%
%TCIMACRO{\tsum }%
%BeginExpansion
{\textstyle\sum}
%EndExpansion
\iota\tau_{j}$ form a basis for the $\mathbb{Z}{}$-module $X^{\ast}(S^{K})$.

Let $B$ be a simple abelian variety over $\mathbb{C}$ of CM-type $(E,\phi)$.
For each $\rho\colon E\rightarrow\mathbb{Q}{}^{\text{al}}$ and $\tau
\in\Gal(\mathbb{Q}{}^{\text{cm}}/\mathbb{Q}{})$, define $\psi_{\rho}%
(\tau)=\phi(\tau^{-1}\circ\rho)$. Then, as $\rho$ runs over the embeddings of
$E$ into $\mathbb{Q}{}^{\text{al}}$, $\psi_{\rho}$ runs over a $\Gamma$-orbit
of CM-types on $\mathbb{Q}{}^{\text{cm}}$. The map $B\mapsto\{\psi_{\rho}\}$
defines a bijection from the set of isogeny classes of simple abelian
varieties of CM-type over $\mathbb{Q}{}^{\text{al}}$ to the set of $\Gamma
$-orbits of CM-types on $\mathbb{Q}{}^{\text{cm}}$. A $\tau\in\Gamma$ fixes
the reflex field $K$ of $B$ if and only if $\tau\phi=\phi.$ Each $\psi_{\rho}$
is the extension to $\mathbb{Q}{}^{\text{cm}}$ of a primitive CM-type on $K$,
and the kernel of $S^{K}\rightarrow\MT(B)$ is the intersection of the kernels
of the $\psi_{\rho}$ (see Milne 1999, \S 2). For any automorphism $\sigma$ of
$\mathbb{C}{}$, $\sigma B$ has reflex field $\sigma K$.

Choose a finite set of $B$'s with reflex field contained in $K$ such that the
corresponding $\psi$'s generate $S^{K}$, and let $A$ be their product. The
canonical map $S^{K}\rightarrow\MT(A)$ is then an isomorphism. According to
(\ref{cm3}), we may take $A$ to be defined over a subfield $k$ of $\mathbb{C}$
of finite degree over $\mathbb{Q}{}$ . Then, if $K$ is Galois over
$\mathbb{Q}{}$, $A_{\ast}=_{\text{df}}\Res_{k/\mathbb{Q}{}}A$ has reflex field
contained in $K$, and the canonical map $S^{K}\rightarrow\MT(A_{\ast})$ is an isomorphism.
\end{proof}

\textsc{Proof of Proposition \ref{g6}. }Let $A$ be an abelian variety over
$\mathbb{Q}{}$. Choose $B$ as in (\ref{z2}). After (\ref{z5}), we may assume
$B$ is defined over a number field $F$. Then $B_{\ast}=%
%TCIMACRO{\tprod _{\sigma\in\Sigma_{F/\mathbb{Q}{}}}}%
%BeginExpansion
{\textstyle\prod_{\sigma\in\Sigma_{F/\mathbb{Q}{}}}}
%EndExpansion
\sigma B$ is defined over $\mathbb{Q}{}$. Let $K$ be any CM-subfield of
$\mathbb{Q}{}^{\text{al}}$ containing the reflex field of $(\MT(A\times
B_{\ast}),h^{A\times B_{\ast}})^{\text{ab}}$ and finite and Galois over
$\mathbb{Q}{}$, and choose $C$ as in (\ref{z6}). Consider the abelian variety
$A\times B_{\ast}\times C$. It is defined over $\mathbb{Q}$, and (\ref{z4},
\ref{z1}, \ref{z1c}) show that it is maximal. By assumption, the reflex field
of $(\MT(A\times B_{\ast}\times C),h^{A\times B_{\ast}\times C})^{\text{ab}}$
is contained $K$, and so the canonical homomorphism $S\rightarrow\MT(A\times
B_{\ast}\times C)^{\text{ab}}$ factors through $S^{K}$. But $\MT(A\times
B_{\ast}\times C)$ surjects onto $\MT(C)\cong S^{K}$, and so $S^{K}%
\cong\MT(A\times B_{\ast}\times C)^{\text{ab}}$. This completes the proof of
Proposition \ref{g6}.

\subsection{Characterizing $P^{\mathcal{H}{}}\rightarrow P^{\CM}$}

\begin{theorem}
\label{g11}Let $A$ be an abelian variety over $\mathbb{Q}{}$. If
$\MT(A)^{\text{der}}$ is simply connected and $\MT(A)^{\text{ab}}\cong S^{K}$
for some $K$, then the isomorphism class of $P^{A}\rightarrow P^{\CM,K}$ is
uniquely determined by its class over $\mathbb{R}{}$.
\end{theorem}

\begin{proof}
According to Proposition \ref{t4} and (\ref{period3}), the isomorphism classes
in $\mathsf{P}(G_{B}^{A}\rightarrow G_{B}^{\CM,K};P^{\CM,K})$ are classified
by $H^{1}(\mathbb{Q}{},H)$, where $H$ is the twist of $\MT(A)^{\text{der}}$ by
the period torsor. Therefore $H$ is simply connected, and so this follows from
the theorem of Kneser, Harder, and Chernousov (\ref{g02}).
\end{proof}

Let $\mathcal{H}{}$ denote the class of abelian varieties $A$ over
$\mathbb{Q}{}$ such that $\MT(A)^{\text{ad}}$ has no factor of type
$D^{\mathbb{H}{}}$. Note that, because of (\ref{g06}\ref{d1}), any abelian
variety for which $\MT(A)^{\text{der}}$ is simply connected lies in
$\mathcal{H}{}$. The next result shows that the abelian varieties${}$
satisfying the conditions of (\ref{g11}) are cofinal among all abelian
varieties in $\mathcal{H}{}{}$.

\begin{proposition}
[Blasius and Borovoi 1999, 3.2]\label{g16}Let $A$ be an abelian variety in
$\mathcal{H}{}{}$. For any sufficiently large CM-field $K$, there exists an
abelian variety $B$ over $\mathbb{Q}{}$ such that $\MT(A\times B)^{\text{ab}%
}\cong S^{K}$ and $\MT(A\times B)^{\text{der}}$ is simply connected$.$
\end{proposition}

\begin{proof}
When $A$ is in $\mathcal{H}{}$, the proof of Proposition \ref{g6} can be
modified to show that $B$ can be chosen in such a way that $\MT(A\times
B)^{\text{der}}$ is simply connected --- in fact, this significantly
simplifies the proof.
\end{proof}

\begin{theorem}
\label{g17}The isomorphism class of $P^{\mathcal{H}{}}\rightarrow P^{\CM}$ is
uniquely determined by its class over $\mathbb{R}{}$.
\end{theorem}

\begin{proof}
Similar to that of Theorem \ref{g7}.
\end{proof}

\begin{note}
\label{g24}The main theorem of Blasius and Borovoi 1999 (Theorem 1.5) states
that the isomorphism class of $P^{\mathcal{H}{}}\rightarrow P^{\CM}$ is
determined by the cohomology class of $P^{\mathcal{H}{}}$ in $\varprojlim
_{A\in\mathcal{H}{}}H^{1}(\mathbb{R}{},G_{\text{dR}}^{A})$. As noted in the
introduction, this statement is false.\footnote{They prove (ibid. 5.2) that
their conditions determine the cohomology class of $P^{\mathcal{H}{}}$ in
$\varprojlim_{A\in\mathcal{H}{}}H^{1}(\mathbb{Q}{},(G_{\text{dR}}^{A})^{\circ
})$, whereas what is needed is that the conditions determine the cohomology
class of $P^{\mathcal{H}{}}\rightarrow P^{\CM}$ in $H^{1}(\mathbb{Q}%
{},(G_{\text{dR}}^{\mathcal{H}{}})^{\circ\text{der}})$. There are
homomorphisms%
\[
H^{1}(\mathbb{Q}{},(G_{\text{dR}}^{\mathcal{H}{}})^{\circ\text{der}%
})\rightarrow H^{1}(\mathbb{Q}{},(G_{\text{dR}}^{\mathcal{H}{}})^{\circ
})\rightarrow\varprojlim_{A\in\mathcal{H}{}}H^{1}(\mathbb{Q}{},(G_{\text{dR}%
}^{A})^{\circ}),
\]
but they are not injective.} Similarly, their statement (ibid 1.6) concerning
the existence of an analogue of their Theorem 1.5 for the category of all
motives is false. The correct statement is Theorem \ref{c8} of this paper.
\end{note}

\subsection{Characterizing $P^{\AV}\rightarrow P^{\CM}$ in terms of a lifting
property}

Let $A$ be an variety over $\mathbb{Q}{}$ such that $\MT(A)\cong S^{K}$ for
some CM-field $K$. Recall that $\MT(A)^{\text{der}}=\Ker(G^{A{}}\rightarrow
G_{\text{B}}^{\CM,K})$. We say that $P^{A}\rightarrow P^{\CM,K}$ has the
\emph{lifting property }if there exists a surjective homomorphism $\tilde
{G}\rightarrow G_{\text{B}}^{A}$such that the kernel of $\tilde{G}\rightarrow
G_{\text{B}}^{\CM,K}$ is the universal covering group of $\MT(A)^{\text{der}}%
$and $P^{A}$ lifts to $\tilde{G}$.

\begin{theorem}
\label{g18}Let $A$ be an abelian variety over $\mathbb{Q}{}$ such that
$\MT(A)^{\text{der}}$ satisfies (*)${}$ and $S^{K}\cong\MT(A)^{\text{ab}}$ for
some CM-subfield $K$ of $\mathbb{C}$. Up to isomorphism, there exists at most
one object $P\rightarrow P^{\CM,K}$ in $\mathsf{P}(G^{A}\rightarrow
G_{\text{B}}^{\CM,K};P^{\CM,K})$ such that

\begin{enumerate}
\item $P\rightarrow P^{\CM,K}$ has the lifting property, and

\item $(P\rightarrow P^{\CM,K})_{\mathbb{R}{}}\approx(P^{A}\rightarrow
P^{\CM,K})_{\mathbb{R}{}}$.
\end{enumerate}
\end{theorem}

\begin{proof}
Apply Proposition \ref{g04} to $H=\MT(A)^{\text{der}}$.
\end{proof}

The notion of a lifting property extends in an obvious fashion to infinite
sets of abelian varieties.

\begin{theorem}
\label{g13}Up to isomorphism, there exists at most one object $P\rightarrow
P^{\CM}$ in $\mathsf{P}(G_{\text{B}}^{\AV}\rightarrow G_{\text{B}}%
^{\CM};P^{\CM})$ such that

\begin{enumerate}
\item $P\rightarrow P^{\CM}$ has the lifting property, and

\item $(P\rightarrow P^{\CM})_{\mathbb{R}{}}\approx(P^{\AV}\rightarrow
P^{\CM})_{\mathbb{R}{}}$.
\end{enumerate}
\end{theorem}

\begin{proof}
Similar to the proof of Theorem \ref{g7}.
\end{proof}

\begin{remark}
\label{g19}Deligne's hope that all Shimura varieties with rational weight are
moduli varieties for motives (Deligne 1979a, p248) implies that $P^{\AV}%
\rightarrow P^{\CM}$ has the lifting property over $\mathbb{C}{}$. It would be
interesting to prove this unconditionally.
\end{remark}

\subsection{Characterizing $P^{\AV}\rightarrow P^{\Art}$}

\begin{theorem}
\label{g20}Let $A$ be an abelian variety over $\mathbb{Q}{}$ such that
$\MT(A)^{\text{ab}}\cong S^{K}$ for some CM-subfield $K$ of $\mathbb{C}{}$.

\begin{enumerate}
\item If $\MT(A)^{\text{der}}$ is simply connected, then the isomorphism class
of $P^{A}\rightarrow P^{\Art}$ is uniquely determined by its class over
$\mathbb{R}{}$.

\item If $\MT(A)^{\text{der}}$ satisfies (*), then, up to isomorphism, there
exists at most one object $P\rightarrow P^{\Art}$ in $\mathsf{P}(G_{\text{B}%
}^{A}\rightarrow\varGamma;P^{\Art})$ having the lifting property and
isomorphic to $P^{A}\rightarrow P^{\Art}$ over $\mathbb{R}{}$.
\end{enumerate}
\end{theorem}

\begin{proof}
(a) According to Proposition \ref{t4}, the isomorphism classes in
$\mathsf{P}(G_{\text{B}}^{A}\rightarrow\varGamma;P^{\Art})$ are classified by
$H^{1}(\mathbb{Q}{},G)$, where $G$ is the twist of
\[
\MT(A)\overset{\ref{period1}}{=}\Ker(G_{\text{B}}^{A}\rightarrow\varGamma)
\]
by the period torsor. Then $G^{\text{ab}}={}_{f}S^{K}$, and so $G$ satisfies
the hypotheses of (\ref{g23}b).

(b) In this case, $G$ satisfies the hypotheses of (\ref{g05}).
\end{proof}

\begin{remark}
\label{g12}For any abelian variety $A$ over $\mathbb{Q}{}$, the map
$G_{\text{dR}}^{A}(\mathbb{R}{})\rightarrow{}_{f}\varGamma(\mathbb{R}{})$ is
surjective (Blasius and Borovoi 1999, 4.4), and so
\[
H^{1}(\mathbb{R}{},(G_{\text{dR}}^{A})^{\circ})\rightarrow H^{1}(\mathbb{R}%
{},G_{\text{dR}}^{A})
\]
is injective. Therefore, for $P\rightarrow P^{\Art}$ in $\mathsf{P}%
(G_{\text{B}}^{A}\rightarrow\varGamma;P^{\Art})$, the isomorphism class of
$P\rightarrow P^{\Art}$ over $\mathbb{R}{}$ is determined by the isomorphism
class of $P$ over $\mathbb{R}{}$.
\end{remark}

\begin{theorem}
\label{g21}(a) There are uncountably many isomorphism classes in
$\mathsf{P}(G_{\text{B}}^{\mathcal{H}{}}\rightarrow\varGamma;P^{\Art})$ that
become isomorphic to $P^{\mathcal{H}{}}\rightarrow P^{\Art}$ over
$\mathbb{R}{}$.

(b) If there exists one object $P\rightarrow P^{\Art}$ in $\mathsf{P}%
(G_{\text{B}}^{\AV}\rightarrow\varGamma;P^{\Art})$ having the lifting property
and isomorphic to $P^{\AV}\rightarrow P^{\Art}$ over $\mathbb{R}{}$, then
there are uncountably many.
\end{theorem}

\begin{proof}
(a) According to Proposition \ref{t4}, the isomorphism classes in
$\mathsf{P}(G_{\text{B}}^{\mathcal{H}{}}\rightarrow\varGamma;P^{\Art})$ are
classified by $H^{1}(\mathbb{Q}{},(G_{\text{dR}}^{\mathcal{H}})^{\circ})$,
where $(G_{\text{dR}}^{\mathcal{H}})^{\circ}$ is the twist of
\[
(G_{\text{B}}^{\mathcal{H}})^{\circ}\overset{\ref{period1}}{=}\Ker(G_{\text{B}%
}^{A}\rightarrow\varGamma)
\]
by the period torsor. For an $A\in\mathcal{H}{}$ with $\MT(A)^{\text{ab}}\cong
S^{K}$, consider%
\[
0\rightarrow(G_{\text{dR}}^{A})^{\circ\text{der}}(\mathbb{Q}{})\rightarrow
(G_{\text{dR}}^{A})^{\circ}(\mathbb{Q}{})\rightarrow{}_{f}S^{K}(\mathbb{Q}%
{})\rightarrow H^{1}(\mathbb{Q}{},(G_{\text{dR}}^{A})^{\circ\text{der}}).
\]
Because $(G_{\text{dR}}^{\mathcal{H}{}})^{\circ\text{der}}$ is a countable
product of algebraic groups,%
\[
\varprojlim\nolimits^{1}(G_{\text{dR}}^{A})^{\circ\text{der}}(\mathbb{Q}%
{})=0=\varprojlim\nolimits^{1}H^{1}(\mathbb{Q}{},(G_{\text{dR}}^{A}%
)^{\circ\text{der}})
\]
(by \ref{c3}). Therefore,%
\[
\varprojlim\nolimits^{1}(G_{\text{dR}}^{A})^{\circ}(\mathbb{Q}{}%
)\cong\varprojlim\nolimits^{1}{}_{f}S^{K}(\mathbb{Q}{}),
\]
which is uncountable (\ref{p8}). Similarly,%
\[
\varprojlim\nolimits^{1}(G_{\text{dR}}^{A})^{\circ}(\mathbb{R}{}%
)\cong\varprojlim\nolimits^{1}{}_{f}S^{K}(\mathbb{R}{}),
\]
which is zero because $_{f}S_{/\mathbb{R}{}}^{K}$ is a product of copies of
$\mathbb{G}_{m}$. In view of Proposition \ref{c1}, this completes the proof.

(b) Similar.
\end{proof}

\begin{remark}
\label{g22}I expect that there are uncountably many distinct isomorphism
classes in $\mathsf{P}(G_{\text{B}}^{\AV}\rightarrow\varGamma;P^{\Art})$ that
become equal to the class of $P^{\AV}\rightarrow P^{\Art}$ over every field
$\mathbb{Q}{}_{l}$.
\end{remark}

\subsection{Characterizing $P^{\AV}$}

\begin{theorem}
\label{mt}Let $A$ be an abelian variety over $\mathbb{Q}{}$ such that
$\MT(A)^{\text{ab}}\cong S^{K}$ for some CM-field $K$.

\begin{enumerate}
\item If $\MT(A)^{\text{der}}$ is simply connected, then, up to isomorphism,
$P^{A}$ is the only $G_{\text{B}}^{A}$-torsor $P$ such that

\begin{enumerate}
\item $P\wedge^{G_{\text{B}}^{A}}\varGamma\approx P^{\Art}$, and

\item $P_{/\mathbb{R}{}}\approx P_{/\mathbb{R}{}}^{A}$.
\end{enumerate}

\item If $\MT(A)^{\text{der}}$ satisfies (*), then, up to isomorphism, there
exists at most one $G_{\text{B}}^{A}$-torsor $P$ such that

\begin{enumerate}
\item $P\wedge^{G_{\text{B}}^{A}}\varGamma\approx P^{\Art}$,

\item there exists a surjective homomorphism $u\colon G^{\prime}\rightarrow
G_{\text{B}}^{A}$ such that $u^{\text{der}}\colon G^{\prime\text{der}%
}\rightarrow(G_{\text{B}}^{A})^{\text{der}}$ is the universal covering group
of $(G_{\text{B}}^{A})^{\text{der}}$ and $P=uP^{\prime}$ for some $G^{\prime}%
$-torsor, and

\item $P_{/\mathbb{R}{}}\approx P_{/\mathbb{R}{}}^{A}$.
\end{enumerate}
\end{enumerate}
\end{theorem}

\begin{proof}
(a) Let $P$ satisfy the conditions (i) and (ii), and choose a morphism
$P\rightarrow P^{\Art}$. Then Theorem \ref{g20} and Remark \ref{g12} show that
$(P\rightarrow P^{\Art})\approx(P^{A}\rightarrow P^{\Art})$.

(b) Similar to (a).
\end{proof}

\subsection{Other fields}

Hasse principles are known to hold for some fields other than number fields.
For example, Scheiderer (1996) proves a Hasse principle for the Galois
cohomology groups of connected linear algebraic groups over perfect fields
with virtual cohomological dimension $\leq1$. However, if the field $k$ is not
countable, the affine group scheme $G$ attached to the category of abelian
motives over $k$ will not be a countable inverse limit of algebraic groups. In
particular, the relation of the flat cohomology group of $G$ to the Galois
cohomology groups of its algebraic quotients is unknown, and so (\emph{pace
}Flicker 2001) such results do not imply Hasse principles for period torsors.

\newpage

%\newpage
\pagestyle{myheadings}\markboth{}{} \setlength{\textwidth}{5.5in}
\setlength{\oddsidemargin}{.5in} \setlength{\evensidemargin}{.5in}
\setlength{\parindent}{-.25in} \setlength{\textheight}{10in}

\addcontentsline{toc}{section}{References}

\section*{References}

{\small \vspace{-0.1in}\hfill}

{\small \textsc{Atiyah, M. F. }1961.\textsf{ }Characters and cohomology of
finite groups. Inst. Hautes \'{E}tudes Sci. Publ. Math. No. 9, 23--64. }

{\small \textsc{Blasius, D., and M. Borovoi. }1999.\textsc{ }On period
torsors. In Automorphic forms, automorphic representations, and arithmetic
(Fort Worth, TX, 1996), 1--8, Proc. Sympos. Pure Math., 66, Part 1, Amer.
Math. Soc., Providence, RI. }

{\small \textsc{Bousfield, A. K., and D. M. Kan.} 1972. Homotopy limits,
completions and localizations. Lecture Notes in Mathematics, Vol. 304.
Springer-Verlag, Berlin-New York. }

{\small \textsc{Breen, L. }1990.\textsc{ }Bitorseurs et cohomologie non
ab\'{e}lienne. In The Grothendieck Festschrift, Vol. I, 401--476, Progr.
Math., 86, Birkh\"{a}user Boston, Boston, MA. }

{\small \textsc{Bucur, I., and A. Deleanu.} 1968. Introduction to the theory
of categories and functors. With the collaboration of Peter J. Hilton and
Nicolae Popescu. Pure and Applied Mathematics, Vol. XIX Interscience
Publication John Wiley \& Sons, Ltd., London-New York-Sydney. }

{\small \textsc{Deligne, P}. 1971. Travaux de Shimura. S\'{e}minaire Bourbaki,
23\`{e}me ann\'{e}e (1970/71), Exp. No. 389 (Preliminary version). }

{\small \textsc{Deligne, P. }1979a. Vari\'{e}t\'{e}s de Shimura:
interpr\'{e}tation modulaire, et techniques de construction de mod\`{e}les
canoniques. In Automorphic forms, representations and $L$-functions (Proc.
Sympos. Pure Math., Oregon State Univ., Corvallis, Ore., 1977), Part 2, pp.
247--289, Proc. Sympos. Pure Math., XXXIII, Amer. Math. Soc., Providence,
R.I.. }

{\small \textsc{Deligne, P.} 1979b. Valeurs de fonctions $L$ et p\'{e}riodes
d'int\'{e}grales. With an appendix by N. Koblitz and A. Ogus. Proc. Sympos.
Pure Math., XXXIII, Automorphic forms, representations and $L$-functions
(Proc. Sympos. Pure Math., Oregon State Univ., Corvallis, Ore., 1977), Part 2,
pp. 313--346, Amer. Math. Soc., Providence, R.I., 1979. }

{\small \textsc{Deligne, P}. 1982a. (Notes by J. Milne), Hodge cycles on
abelian varieties. In Hodge cycles, motives, and Shimura varieties. Lecture
Notes in Mathematics, 900. Springer-Verlag, Berlin-New York, pp9--100. }

{\small \textsc{Deligne, P.,} 1982b. Motifs et groupe de Taniyama. In Hodge
cycles, motives, and Shimura varieties. Lecture Notes in Mathematics, 900.
Springer-Verlag, Berlin-New York, pp261--279. }

{\small \textsc{Deligne, P., and J.S. Milne. }1982. Tannakian categories. In
Hodge cycles, motives, and Shimura varieties. Lecture Notes in Mathematics,
900. Springer-Verlag, Berlin-New York, pp101--228.}

{\small \textsc{Flicker, Y. Z.} 2001. Motivic torsors. Israel J. Math. 122,
61--77.}

{\small \textsc{Giraud, J.} 1971. Cohomologie non ab\'{e}lienne. Die
Grundlehren der mathematischen Wissenschaften, Band 179. Springer-Verlag,
Berlin-New York. }

{\small \textsc{Gray, B. I. }1966. Spaces of the same $n$-type, for all $n$.
Topology 5, 241--243. }

{\small \textsc{Grothendieck, A.} 1966 . On the de Rham cohomology of
algebraic varieties. Inst. Hautes \'{E}tudes Sci. Publ. Math. No. 29, 95--103.
}

{\small \textsc{Milne, J. S.} 1988. Automorphic vector bundles on connected
Shimura varieties. Invent. Math. 92, no. 1, 91--128. }

{\small \textsc{Milne, J. S.} 1990. Canonical models of (mixed) Shimura
varieties and automorphic vector bundles. Automorphic forms, Shimura
varieties, and $L$-functions, Vol. I (Ann Arbor, MI, 1988), 283--414,
Perspect. Math., 10, Academic Press, Boston, MA. }

{\small \textsc{Milne, J. S. }1999. Lefschetz motives and the Tate conjecture.
Compositio Math. 117, no. 1, 45--76. }

{\small \textsc{Platonov, V. and A. Rapinchuk.} 1994. Algebraic groups and
number theory. Translated from the 1991 Russian original by Rachel Rowen. Pure
and Applied Mathematics, 139. Academic Press, Inc., Boston, MA. }

{\small \textsc{Scheiderer, C.} 1996. Hasse principles and approximation
theorems for homogeneous spaces over fields of virtual cohomological dimension
one. Invent. Math. 125, no. 2, 307--365.}

{\small \textsc{Serre, J.-P. }1964. Cohomologie galoisienne. Cours au
Coll\`{e}ge de France, 1962-1963. Seconde \'{e}dition. Lecture Notes in
Mathematics 5 Springer-Verlag, Berlin-Heidelberg-New York. }

{\small \textsc{Serre, J.-P. }1992. Topics in Galois theory. Lecture notes
prepared by Henri Damon [Henri Darmon]. With a foreword by Darmon and the
author. Research Notes in Mathematics, 1. Jones and Bartlett Publishers,
Boston, MA. }

{\small \textsc{Serre, J.-P. }1994. Propri\'{e}t\'{e}s conjecturales des
groupes de Galois motiviques et des repr\'{e}sentations $l$-adiques. In
Motives (Seattle, WA, 1991), 377--400, Proc. Sympos. Pure Math., 55, Part 1,
Amer. Math. Soc., Providence, RI. }

{\small \textsc{Tate, J. }1976. Relations between $K_{2}$ and Galois
cohomology. Invent. Math. 36, 257--274. }

{\small \textsc{Wei, W. }1994. \hfill Moduli fields of CM-motives applied to
Hilbert's 12-th problem, \hfill Preprint
http://www.mathematik.uni-bielefeld.de/sfb343/preprints/pr94070.ps.gz. }

{\small \textsc{Wintenberger, J.-P.} 1990. Torseurs pour les motifs et pour
les repr\'{e}sentations $p$-adiques potentiellement de type CM, Math. Ann.
288, no.~1, 1--8. }

\vspace{.3in}

{\small 2679 Bedford Rd, Ann Arbor, MI 48104, USA. }

{\small E-mail: math@jmilne.org; }

{\small Home page: www.jmilne.org/math/ }

\end{document}

%% file: def01.tex
\newcounter{junk}

\newenvironment{mylist}
{\begin{list}
{--}
{\setlength{\leftmargin}{.5in}\setlength{\rightmargin}{.5in}}}
{\end{list}}

\renewenvironment{itemize}
{\begin{mylist}}
{\end{mylist}}

\def\1{{1\mkern-7mu1}}  %% since alas there's no "\Bbb{1}";
\newcommand\ad{\textrm{ad}}
\newcommand\Ad{\operatorname{Ad}}

\newcommand\Aut{\operatorname{Aut}}

\newcommand\Gal{\operatorname{Gal}}
\newcommand\GL{\operatorname{GL}}

\newcommand\Hom{\operatorname{Hom}}

\newcommand\Ker{\operatorname{Ker}}

\newcommand\MT{\operatorname{MT}}
\newcommand\Nm{\operatorname{Nm}}

\newcommand\ord{\operatorname{ord}}

\newcommand\Res{\operatorname{Res}}
\newcommand\SO{\operatorname{SO}}
\newcommand\Spec{\operatorname{Spec}}

\def\AV{{\mathsf{AV}}}%abelian varieties
\def\Art{{\mathsf{Art}}}%Artin motives
\def\CM{{\mathsf{CM}}}%CM-motives
\def\Mot{{\mathsf{Mot}}}%motives
\def\ad{{\mathrm{ad}}}

\def\MT{{\mathrm{MT}}}

\def\ord{{\mathrm{ord}}}

\let\xr=\xrightarrow